\pdfoutput=1
\newif\ifelsarticle
\elsarticlefalse 
\newif\ifsup
\suptrue 
\ifelsarticle
\documentclass[preprint,12pt]{elsarticle}

\journal{Linear Algebra and its Applications}

\else
  \documentclass[12pt]{article}
  \usepackage{authblk} 
\fi
\usepackage{ifthen}
\usepackage{graphicx}
\usepackage{epsfig} 
\usepackage{epstopdf}
\usepackage{times} 
\usepackage[utf8]{inputenc}
\usepackage[english]{babel}
\usepackage{fancyhdr}
\usepackage{dsfont}
\usepackage{algorithm}
\usepackage{booktabs}
\usepackage{algpseudocode}
\usepackage{changepage}
\usepackage{multirow}

\newboolean{longver}
\setboolean{longver}{true}

\usepackage{amsmath}
\usepackage[table,xcdraw,dvipsnames]{xcolor}

\usepackage{amssymb}

\usepackage{amsthm}

\usepackage{enumitem}
\usepackage{tikz}
\usepackage{pgfplots}
\usetikzlibrary{fit}
\usetikzlibrary{shapes,arrows}
\usetikzlibrary{matrix,decorations.pathreplacing,arrows.meta}
\usepackage{lipsum}
\usepackage[caption=false, font=normalsize, labelfont=scriptsize, textfont=scriptsize]{subfig}

\pgfplotsset{compat=newest}
\usepackage{pgfplotstable}
\usepgfplotslibrary{colormaps}

\usepackage{color}

\newtheorem{thm}{Theorem}

\newtheorem{prop}{Proposition}
\newtheorem{lem}{Lemma}
\newtheorem{rem}{Remark}
\newtheorem{cor}{Corollary}

\newtheorem{cex}{Counterexample}

\newtheorem{ex}{Example}

\newcommand{\R}{\mathbb{R}}

\renewcommand{\S}{\mathbb{S}}
\newcommand{\N}{\mathbb{N}}

\newcommand{\eps}{\varepsilon}

\newcommand{\rank}{\mathop{\mathrm{rank}}}
\newcommand{\diag}{\mathop{\mathrm{diag}}}

\DeclareMathOperator{\Tr}{Tr}

\usepackage{stackengine}
\newcommand\marksymb[2]{\,\stackon{#1}{\scriptscriptstyle (#2)}\, }
\stackMath
\begin{document}
\ifelsarticle
\begin{frontmatter}



\title{
Approximate Projections onto the Positive Semidefinite Cone Using Randomization
}
 \author[label1]{Morgan Jones}
 \affiliation[label1]{organization={School of Electrical and Electronic Engineering},
             addressline={The University of Sheffield},
         email={ morgan.jones@sheffield.ac.uk}}
     
      \author[label2]{James Anderson}
     \affiliation[label2]{organization={Department of Electrical Engineering},
     	addressline={Columbia University},
     	email={ anderson@ee.columbia.edu}}




\else
  \title{Approximate Projection onto the Positive Semidefinite Cone Using Randomization}
  \date{}
\author{Morgan Jones
	\thanks{ \footnotesize M. Jones is with the School of Electrical and Electronic Engineering,
		The University of Sheffield. E-mail: {\tt \small morgan.jones@sheffield.ac.uk} }
	\quad James Anderson%
	\thanks{\footnotesize J. Anderson is with  Department of Electrical Engineering, Columbia University . E-mail: {\tt \small anderson@ee.columbia.edu} }
}
  \maketitle
  \fi

\begin{abstract}
	This paper presents two algorithms that compute approximate Positive Semidefinite (PSD) projections of real symmetric matrices using Randomized Numerical Linear Algebra (RNLA). Classical PSD projection of an $n\times n$ matrix  relies on a deterministic eigen-decomposition with computation that scales as $\mathcal{O}(n^3)$. Our approach leverages RNLA to construct low-rank matrix approximations before projection, significantly reducing the required numerical resources to $\mathcal{O}(k n^2)$, for some user defined fixed parameter $k$. The first algorithm utilizes random sampling to generate a low-rank approximation, followed by a standard eigen-decomposition on this smaller matrix. The second algorithm enhances this process by introducing a scaling approach that aligns the leading-order singular values with the positive eigenvalues, biasing the low-rank approximation to focus on capturing the essential information about the positive eigenvalues for PSD projection. Both methods offer a trade-off between accuracy and computational speed, supported by probabilistic error bounds. {Numerical experiments on large-scale matrices ( $n\approx 20K$) demonstrate that the proposed randomized algorithms effectively approximate PSD projections.}
	
	
\end{abstract}

\ifelsarticle
\begin{highlights}
\item Two methods for projecting onto the positive semidefinite cone of matrices powered by randomized numerical linear algebra.
\item Probabilistic error bounds in both the spectral and Frobenius norms for our proposed projection methods.
\item Numerical examples validate our error bounds and illustrate the advantages and disadvantages of the two proposed projection methods. 
\end{highlights}


\end{frontmatter}
\else
\fi



 \section{Introduction}

 {Many applications in scientific computing theoretically require matrices to be Positive Semidefinite (PSD). Examples include \emph{semidefinite programming (SDP)}~\cite{yurtsever2021scalable}, statistical inference (e.g., estimating correlation matrices)~\cite{higham2002computing}, solutions to Lyapunov equations~\cite{banno2023data}, ellipsoidal shape matrices~\cite{halder2018parameterized} and Hessians of convex functions~\cite{thacker1989role}. However, when these matrices are estimated from noisy data~\cite{higham2002computing} or computed numerically with floating point errors they frequently lose this property and become ``almost'' PSD. Repairing these matrices requires us to carry out} the task of finding the closest PSD matrix, or equivalently projecting a matrix onto the PSD cone. While projecting a matrix onto the PSD cone is a fundamental computation in numerical linear algebra, analyzing the error of such projections is non-trivial~\cite{goulart2020accuracy}. Unfortunately, existing algorithms for computing such projections scale poorly with the matrix dimension, making them impractical for large-scale problems.  To address this limitation, in this paper we provide methods that trade off accuracy for efficiency by employing techniques from Randomized Numerical Linear Algebra (RNLA) to compute approximate PSD projections.

 The problem of projecting a real symmetric matrix $X \in \mathbb{R}^{n \times n}$ onto the PSD cone can be formulated as the convex optimization problem:  $\min_{Y \succcurlyeq 0} \|X-Y\|^2_2$, where $\|\cdot\|_2$ is the Frobenius norm. This problem has a closed form solution: compute the eigen-decomposition of $X$,  zero out the negative eigenvalues, and re-multiply the modified factors to yield the projected matrix. For a dense $n \times n$ matrix, this demands $\mathcal{O}(n^3)$ Floating-Point Operations (FLOPs)~{\cite{pan1999complexity}}, making it impractical for large $n$.  Existing approximate PSD projection methods, such as~\cite{francisco2017fixed} and~\cite{higham1988computing}, also entail $\mathcal{O}(n^3)$ FLOPs. In this work, we are concerned with use cases where $n$ is sufficiently large that a cubic dependence is intractable. We propose the use of RNLA~\cite{halko2011finding} for computing approximate low-rank projections onto the PSD cone. 
 
 Our method of using low rank approximations to construct approximate PSD projections is justified {by the observation that matrices arising in many problem settings are well approximated by low-rank matrices, with rapidly decaying singular values and correspondingly small truncation error. This is evidenced } by the pervasiveness of {approximate} low rank matrices appearing in practical applications and data science as a whole~\cite{udell2019big}. {Moreover, for SDPs with a relatively small number of constraints, it is known that there always exists a low-rank optimal solution~\cite{majumdar2020recent,pataki1998rank,barvinok1995problems}.} Using RNLA we demonstrate the possibility of computing inexact PSD projections with {quadratic rather than cubic complexity. More specifically, our PSD projection method scales with} $\mathcal{O}(kn^2)$ FLOPs, where {$k \ll n$} is a user-defined parameter that can be increased to improve accuracy. At little extra computational cost, we introduce a scaling method that biases the RNLA algorithm to construct a low-rank matrix approximation focused on capturing information about the positive eigenvalues, rather than the leading singular values. We demonstrate that in some cases, this approach yields a better PSD projection approximation. This improvement arises because the method reduces computational effort spent on capturing information on large negative eigenvalues, which are inevitably set to zero in the PSD projection.

 {The PSD projection step is the primary computational bottleneck in many algorithms, particularly first-order methods for solving SDPs, such as those of~\cite{o2016conic, yang2015sdpnal+, zheng2020chordal, zheng2017exploiting}. The cost of PSD projection dominates the per-iteration time complexity of first order SDP solvers~\cite{rontsis2022efficient}. While second-order interior-point methods offer superior convergence rates and accuracy, they are even less scalable due to the need to compute, store, and solve large-scale Hessian systems~\cite{majumdar2020recent}. These limitations highlights the need for scalable approximate PSD projection techniques that significantly reduce computational cost while retaining sufficient accuracy for implementation purposes. In this work, we focus exclusively on the PSD projection subroutine itself. Our proposed method may be used as a drop-in replacement for the exact PSD projection step within existing first-order SDP solvers; however, in this paper we do not study end-to-end SDP convergence and performance, instead concentrating on the efficiency and accuracy of our proposed approximate projections in isolation.}
 
 {Randomization has been employed in prior work to manage the computational challenges of SDPs. For instance, \cite{yuan2018random} proposes an alternating projection scheme where each iteration projects onto a randomly selected subset of the constraints. Similarly, \cite{guedes2024sparse} uses random sketching to construct lower-dimensional surrogate SDPs. In contrast, our use of randomization is independent to SDPs, instead, we employ randomness for PSD cone projection, allowing for the potential of seamless integration into existing first-order SDP solvers.} The works of~\cite{fawzi2021faster, rontsis2022efficient, souto2022exploiting,kang2025factorization} demonstrate the usefulness of using approximate PSD projections when solving SDPs. In~\cite{fawzi2021faster}, an approximate PSD projection based Jacobi's method showed impressive performance but degraded significantly with larger matrices. On the other hand, both~\cite{rontsis2022efficient} and~\cite{souto2022exploiting} base their projection on Krylov methods, which are known to {exhibit} numerical instability issues in finite-precision arithmetic~\cite{halko2011finding}. {In particular, the Lanczos algorithm, widely used for constructing Krylov subspaces, does not fulfil its theoretical guarantees when implemented in finite-precision arithmetic~\cite{meurant2006lanczos}. }Furthermore, Krylov methods need multiple passes over the matrix, a computational bottleneck for large matrices that do not fit in fast memory, while RNLA methods require just a single pass~\cite{halko2011finding}.  {The work of~\cite{kang2025factorization} bases its PSD projection on polynomial approximations of the ReLU formulation, demonstrating impressive results, but ultimately still requires dense matrix multiplication at a cost of $\mathcal{O}(n^3)$ flops.}
 

\section{Notation} 
We denote the set of {real} square symmetric matrices by ${\S_n}$ and the set of {real} square symmetric positive semidefinite matrices by $\S_n^+$. For $A \in \S_n^+$,  $\sqrt{A} \in \S_n^+$ is defined such that $\sqrt{A}^\top \sqrt{A}=\sqrt{A}\sqrt{A}=A$ (note that the square root of a matrix, defined in this way, is unique~\cite{higham1997stable}. We also use  $A^{\frac{1}{2}}$ in place of $\sqrt{A}$). For $A \in \S_n$, we index the eigenvalues of $A$ in increasing order by $\lambda_1(A) \ge \dots \ge \lambda_n(A)$ and the singular values by $\sigma_1(A) \ge \dots \ge \sigma_n(A)$.  We denote the Schatten $p$-norms for $A \in \S_n$ as $\|A\|_p=\left( \sum_{i=1}^n \sigma_i(A)^p \right)^{\frac{1}{p}}$. In particular, the Frobenius norm and spectral norm of $A$ are denoted by $\|A\|_2=\sqrt{\Tr(A^\top A)}=\sqrt{\sum_{i=1}^n \sigma_i(A)^2}$ and $\|A\|_{\infty}=\max_{1 \le i \le n} \sigma_i(A)=\sigma_1(A)$ under our $p$-norm notation. When it is not important which norm is used we simply use $\|\cdot\|$.  We denote $D=\diag(\lambda_1,\dots,\lambda_n)$ to mean that $D \in \R^{n \times n}$ is a diagonal matrix such that $D_{i,i}=\lambda_i$ and $D_{i,j}=0$ for $i \ne j$. Moreover, for diagonal matrices $D=\diag(\lambda_1,\dots,\lambda_n)$ we denote $\max(D,\beta):=\diag(\max\{\lambda_1,\beta\},\dots,\max\{\lambda_n,\beta\})$ and $D_+:=\max(D,0)$. {For a general matrix, $X \in \S_n$, we generalize the $\cdot_+$ operator by defining $X_+$ as the PSD projection given in Eq.~\eqref{eq: analytical PSD proj using eigenvalue}}. We denote the row vector of ones by $\mathbf{1}_{n}=[1,\dots,1] \in \R^{1 \times n}$. {Throughout the paper, we use the hat symbol to indicate that a quantity is an approximation. Specifically, if \( Y = \hat{X} \), then \( Y \approx X \) in some appropriate sense. It is important to note that the hat symbol is not an operator. Consequently, expressions like \(\hat{X}_+\) should always be understood as \((\hat{X})_+\), that is, the PSD projection is applied to the approximated matrix \(\hat{X}\).} {For a real number $a \in \R$ and integer $z \in \N$ such that $z<a<z+1$ we define the rounding of $a \in \R$ to the nearest integer as $[a]:=\begin{cases}
		z \text{ if } a < z+0.5 \\
		z+1 \text{ otherwise.}
\end{cases}$}

\section{Projecting onto the PSD Cone} \label{sec: Projecting onto the Cone of Positive Semidefinite Matrices}



Let us now consider the problem of projecting a given matrix $X \in \S_n$ onto $\S_n^+$, i.e.,  finding the closest element in $\S_n^+$ to  $X \in \S_n$. {This task can be formulated as the} optimization problem:

\begin{align} \label{opt: projecting onto PSD}
	X_+:=\arg \min_{Y \in \S_n^+} \|X-Y\|_2^2.
\end{align}
{The solution to Opt.~\eqref{opt: projecting onto PSD} has a closed form solution given by the theorem below.} 
\begin{thm}[Page 399~\cite{boyd2004convex}] \label{thm: Analytical PSD Projection}
	Consider $X \in \S_n$. Then there is a unique solution {$X_+$ to optimization}~\eqref{opt: projecting onto PSD} given by
	\ifthenelse{\boolean{longver}}{}{\vspace{-0cm}}\begin{align} \label{eq: analytical PSD proj using eigenvalue}
		X_+=U D_+ U^\top,
	\end{align}
	where $X =  UDU^\top $ is the eigen-decomposition of $X \in \S_n$.
\end{thm}
 {Note that Equation~\eqref{eq: analytical PSD proj using eigenvalue} is not only a closed form solution to Optimization Problem~\eqref{opt: projecting onto PSD}, where the objective is taken with respect to the Frobenius norm, denoted in this paper as the Schatten 2-norm $\| \cdot \|_2$, but also the closed form solution with respect to any unitarily invariant norm (see Lemma 2.1~\cite{goulart2020accuracy}).}

\begin{figure}
		\subfloat[ {\scriptsize For a symmetric matrix $X=\begin{bmatrix} x & y \\ y & z\end{bmatrix}$ the PSD cone, $\S_2^+$, can be characterized as the set of matrix element points obeying $z \ge 0$, $x \ge 0$ and $xz \ge y^2$ by ensuring the leading principal minors are nonnegative. 
  } \label{fig: PSD proj}]{		\begin{tikzpicture}[scale=0.8]
		\begin{axis}[
			colorbar style={ ztick={0.1,0.50, 0.75,  2.00} },
			view={310}{40},
			point meta min=0, 
			point meta max=7, 
			xlabel={$x$},
			ylabel={$y$},
			zlabel={$z$},
			xtick={0,0.5,1},
			ytick={-1,0,1},
			zmax=1,
			zmin=0,
			ymin=-1, 
			ymax=1,
			xmax=1,
			xmin=0,
				3d box,
			colormap/hot2,
			domain=0:1,
			y domain=-1:1
			]
			\addplot3[surf,	shader=interp,	opacity=0.4] 
			{(y^2)/x};
					\fill[red] (0,-0.9, 0.25) circle (2pt) node[below left] {$X$};
					\draw[->,thick,blue] (0,-0.9, 0.25) -- (0.25,0, 0.25) ;
					\node[below, rotate=-50] at (0.3,-0.4,0.25) { \textcolor{blue}{PSD projection}};
						\fill[green!50!black] (0.25,0.02, 0.27)circle (2pt) node[above left] {$X_+$};
						\addlegendentry{$\S_2^+$}
									\addplot3[surf,	shader=interp,	opacity=0] 
						{0};
				\addlegendentry{\tikz[baseline] \node[draw, fill=white, above left] {\text{ }}; $\S_2$}
		\end{axis}
	\end{tikzpicture}  }  \hfill 
		\subfloat[ {\scriptsize The non-expansive property of projections onto closed convex  sets in $\R^n$.} \label{subfig: non expansive}]{	\begin{tikzpicture}[>=Stealth]
		
		\draw[fill=gray!20] (0,0) ellipse (1 and 2);
		
		
		\node at (-1.7,2) {$\text{Closed convex}$};
		\node at (-1.5,1.7) {$\text{set}$};
		
		\fill[red] (2, 2) circle (2pt) node[above ] {$x$};
		\fill[red] (1.5, -1.5) circle (2pt) node[below] {$y$};
		
		\draw[->,blue] (2, 2) -- (0.7, 1.4);
		\draw[->,blue] (1.5, -1.5) -- (0.85, -1);
		
		\fill[green!50!black] (0.7, 1.4) circle (2pt) node[left ] {\small $\pi(x)$};
		\fill[green!50!black] (0.85, -1) circle (2pt) node[ left] {\small $\pi(y)$};
		
		\draw[dashed,red] (2, 2) -- (1.5, -1.5);
		\draw[dashed, green!50!black] (0.7, 1.4) -- (0.85, -1);
		
		\node[red,rotate=85] at (1.6, 0.6) {\scriptsize $d(x,y)$};
		\node[green!50!black,rotate=95] at (0.55, 0.2) {\scriptsize $d(\pi(x),\pi(y))$};
		
				\node[green!50!black]  at (-2.4,0.8) (node1) {\footnotesize $d(\pi(x),\pi(y))$};
						\node[red]   at (-2.6,0.4) (node2){\footnotesize $\textcolor{black}{<}d(x,y)$};
						
						\node[draw, rectangle, fit=(node1) (node2)] {};
		
	\end{tikzpicture}}\hfill 
		
		\vspace{-5pt}
		\caption{\footnotesize Projection descriptions and properties.}
\end{figure}
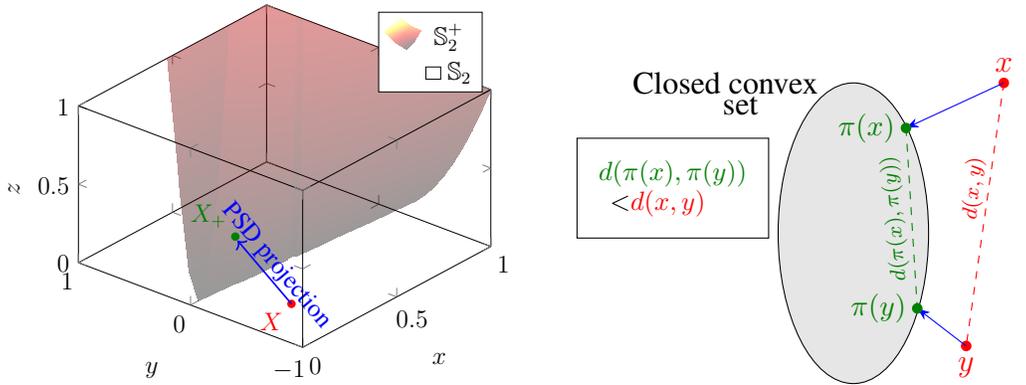

Theorem~\ref{thm: Analytical PSD Projection} (visualized in Fig.~\ref{fig: PSD proj}) shows that finding projections onto the PSD cone by solving Opt.~\eqref{opt: projecting onto PSD} numerically boils down to computing full eigen-decompositions {at a cost of $\mathcal{O}(n^3)$ Floating Point Operations (FLOPs)~\cite{pan1999complexity}.}  We  propose methods that compute approximate projections, outputting $\hat{X}_+ \in \S_n^+$ such that $\hat{X}_+ \approx X_+$, at a cost of $\mathcal{O}(kn^2)$, where $k{\ll n}$ is a user chosen parameter that can be increased to reduce approximation error.

		

		
		
		


\subsection{{\textbf{Low-Rank Randomized PSD Projection}}} 
Techniques from Randomized Numerical Linear Algebra (RNLA) offer an efficient and scalable alternative to classical deterministic numerical linear algebra algorithms. {Surprisingly, RNLA approaches are often more robust to numerical error than their deterministic counterparts, despite the apparent variability introduced by randomization. Moreover, RNLA methods are often easily parallelizable. These properties have been well-documented~\cite{mahoney2011randomized,halko2011finding}.}

{While approximate projection onto the PSD cone was not considered in~\cite{halko2011finding}, we show how to adapt their RNLA-based approach, which approximates the range of a matrix via a lower-dimensional subspace, to derive efficient PSD projection approximations. The methods of~\cite{halko2011finding}} seek to find a  matrix $Q \in \R^{n \times k}$ with orthonormal columns where $Q^\top Q=I$, and $k < n$ such that $X \approx Q Q^\top X$, where $\approx$ should be interpreted as
%
\begin{align} \label{ineq: X approx Q Q X}
	\|(I-QQ^\top)X\|\le \epsilon,
\end{align}
for some user-specified tolerance $\epsilon 
>0$ {(clearly $Q Q^\top \ne I$)}. The fewer columns there are in $Q$ the larger the computational saving. Following~\cite{halko2011finding}, $Q$ can be found via random sampling { (at a  cost proportional to the cost of a matrix-vector multiplication) and an economy QR decomposition of an $n\times k$ matrix.}

{Specifically,  to find $Q\in \R^{n \times k}$ we sample the range of $X$ by drawing an $n$-dimensional random\footnote{ {{The distribution of the test vectors $\omega_i$ is not critical; i.i.d. standard Gaussian vectors will suffice. However, structured matrices like the Subsampled Randomized Fourier Transform (SRFT) are often used in practice~\cite{rokhlin2008fast}. Theoretically, SRFTs improve the computational complexity of the sketching step from $\mathcal{O}(n^2 k)$ to $\mathcal{O}(n^2 \log k)$ operations~\cite{woolfe2008fast}. This speedup comes at the cost of weaker theoretical error bounds: SRFT bounds lack scale-invariant expectation bounds, introduce a dependence on the matrix dimension $n$, require significantly larger oversampling, and exhibit polynomial rather than exponential decay in failure probability (see \cite{halko2011finding}, Sec. 11). Despite these weaker guarantees, SRFTs perform exceptionally well in practice, and our theoretical arguments can be straightforwardly adapted to accommodate them.}}} test vector $\omega_i$ and computing $y_i=X \omega_i$}. Sampling $k$ times and writing in matrix notation gives $Y=X\Omega$ where $\Omega=[\omega_1,\dots, \omega_k]$. The intuition behind randomized algorithms is that the span of the column space of $X$ is approximately the the span of the column space of $Y$, thus an exact basis for the range of $Y$ should provide a good approximation for the basis of $X$. The range of $Y$ is a subspace of dimension $k$ which can be cheaply computed using a standard economic QR decomposition, obtaining $Q\in \R^{n \times k}$. If $X$ has an approximately rank-$k$ structure then Ineq.~\eqref{ineq: X approx Q Q X} holds. In practice, if a rank-$k$ approximation to $X$ is sought, then $k+l$ random sample are obtained {for some oversampling parameter $l \in \N$}. Typically $l=10$ works very well in practice~\cite{halko2011finding}.
It is well known that low-rank approximations perform best when the singular values decay rapidly. For matrices with slow spectral decay, one can instead employ the power iteration method and work with the matrix $W=(XX^\top)^q X$ with $q$ a small integer instead of $X$~{\cite{martinsson2010normalized,halko2011finding}}. Intuitively, $W$ and $X$ have the same singular vectors, however, the singular values of $W$ that are less than one are compressed towards zero. The steps above are formalized as an algorithm referred to as the \emph{randomized range finder} that is recalled from~\cite{halko2011finding} in Alg.~\ref{alg:RSNA}.

 \begin{algorithm}
	\caption{Randomized Range Finder ~\cite[p.244]{halko2011finding}.}\label{alg:RSNA}
	\textbf{Function} $Q=\texttt{range$\_$finder}(X,k,l,q)$. \\
	\hspace*{\algorithmicindent} \textbf{Input:} $X\in \R^{m \times n}$, {approximate rank} $k\in \N$, {oversampling parameter} $l\in \N$\\
	\text{ } \hspace{1.7cm} {and power iteration parameter} $q \in \N$. \\
	\hspace*{\algorithmicindent} \textbf{Output:} $Q\in \R^{n \times (k+l) }$. 
	\begin{algorithmic}[1]
		\State $\Omega = \texttt{randn}(m,k+l)$
		\State $Q=\texttt{orth}( (XX^\top)^q X \Omega)$ {\Comment{\texttt{[Q,$\sim$]=QR($(XX^\top)^q X \Omega$,`econ')}}}
	\end{algorithmic}	
\end{algorithm}

For a symmetric matrix, $X$, we use the approximation $X \approx {\hat{X}:=} QQ^\top X QQ^\top$ rather than $X \approx QQ^\top X$. This form preserves symmetry, allowing us to perform the approximate eigen-decomposition required by our method for computing approximate PSD projections. The approximation $X \approx {\hat{X}:=} QQ^\top X QQ^\top$ can be understood as a simultaneous projection onto subsets of both the row and column spaces of $X$, which are identical for symmetric matrices and hence will have the same projection matrix $Q$. The following result bounds the difference between the original matrix and it's approximate projection ${\hat{X}:=} QQ^\top XQQ^\top$, where $Q$ satisfies Ineq.~\eqref{ineq: X approx Q Q X}. 
\begin{lem} \label{lem: 2 and F bound involving Q}
	Consider $X \in \S_n$ and $Q \in \R^{n \times k}$. Suppose $Q^\top Q=I$, $\|(I-QQ^\top)X\|_\infty \le \eps_1$ and $\|(I-QQ^\top)X\|_2 \le \eps_2$ then
	\begin{align} \label{ineq: 2 norm bound}
		\|X- {\hat{X}}\|_\infty & \le 2\eps_1,\\ \label{ineq: F norm bound}
		\|X- {\hat{X}} \|_2 & \le  {\sqrt{2}}\eps_2,
	\end{align}
{where $\hat{X}:=Q Q^\top X Q Q^\top$}
\end{lem} 
\begin{proof}
	Under any Schatten norm, using the triangle inequality, the sub-multiplication property and the fact that $\|A^\top\|=\|A\|$  {(which implies $\|X(I-QQ^\top)\|=\|(X(I-QQ^\top))^\top\|=\|(I-QQ^\top)X\|$)}, it is easy to show that
	\begin{align}  \nonumber
		&{\| X - \hat{X}\|}=\|X- Q Q^\top X Q Q^\top\| =  \|X-Q Q^\top X + Q Q^\top X -Q Q^\top X Q Q^\top  \| \\ \label{ineq:bound symetrix approx} 
		& \le \|X-Q Q^\top X \|+  \| Q Q^\top( X- X Q Q^\top ) \| \le (1+ \|QQ^\top\|)\|(I-QQ^\top)X\|.
	\end{align}
	Since $Q \in \R^{n \times k}$ it follows $\rank(QQ^\top) \le k$. Moreover, {Since \( QQ^\top \) is a real symmetric matrix, the Spectral Theorem (Page 330~\cite{strang2000linear}) ensures that it has real eigenvalues and an orthonormal basis of eigenvectors: \( QQ^\top v_i= \lambda_i(QQ^\top)v_i \) with \( v_i^\top v_i = 1 \).} Because $Q^\top Q=I$ it follows then that,
	\begin{align} \label{eq: QQ eigenvalues}
		&\left( \lambda_i(QQ^\top) \right)^2=(QQ^\top v_i)^\top QQ^\top v_i=v_i^\top Q Q^\top Q Q^\top v_i\\ \nonumber
		& \qquad =v_i^\top Q  Q^\top v_i=\lambda_i(QQ^\top) \implies \lambda_i(QQ^\top)=0 \text{ or } 1.
	\end{align}
 Hence, 
	\begin{align} \label{eq: QQ 2 norm}
		&\|QQ^\top\|_\infty=\sigma_1(X) = 1.
	\end{align}
	Ineq.~\eqref{ineq: 2 norm bound} follows from applying Ineq.~\eqref{ineq:bound symetrix approx} and Eq.~\eqref{eq: QQ 2 norm}.  
	
	We next prove Ineq.~\eqref{ineq: F norm bound}.  {Using the facts 
	\begin{align} \label{pfeq: facts 1}
	&Q^\top Q=I\\ \label{pfeq: facts 2} 
	&(I-QQ^\top)(I-QQ^\top)=(I-QQ^\top)\\
	& \label{pfeq: facts 3} QQ^\top(I-QQ^\top)=(I-QQ^\top)QQ^\top=0,
	\end{align} and the linearity of the trace operator we get that
	\begin{align} \nonumber
	&	\|X-\hat{X}\|_2^2=\|X-Q Q^\top X + Q Q^\top X -Q Q^\top X Q Q^\top  \|_2^2\\ \nonumber
		&=\|(I-Q Q^\top) X + Q Q^\top X(I- Q Q^\top)  \|_2^2\\ \nonumber
		&=\Tr \hspace{-0.1cm} \bigg( \hspace{-0.15cm}  \left(X(I \hspace{-0.05cm} - \hspace{-0.05cm} Q Q^\top)\hspace{-0.075cm}+\hspace{-0.075cm}(I \hspace{-0.05cm} -\hspace{-0.05cm}  Q Q^\top)XQ Q^\top \right)  \left((I \hspace{-0.05cm} - \hspace{-0.05cm} Q Q^\top) X \hspace{-0.075cm}+ \hspace{-0.075cm} Q Q^\top X(I \hspace{-0.05cm} -\hspace{-0.05cm}  Q Q^\top) \right) \hspace{-0.15cm}  \bigg)\\ \nonumber
		&\marksymb{=}{\eqref{pfeq: facts 1}-\eqref{pfeq: facts 3}}  \Tr(X(I-QQ^\top)X)+\Tr((I-QQ^\top) X QQ^\top X(I-QQ^\top))\\ \label{pfeq: square frob}
		&=\|(I-QQ^\top)X\|_2^2 + \|QQ^\top X(I-QQ^\top)\|_2^2.
	\end{align}
	We next bound the second term in Eq.~\eqref{pfeq: square frob}. Using the linearity of the trace operator and positivity of the norm $-\| \cdot \| \le 0$ we get that,
	\begin{align} \nonumber
&	\|QQ^\top X(I-QQ^\top)\|_2^2=\| X(I-QQ^\top) - (I-QQ^\top)X(I-QQ^\top)\|_2^2\\ \nonumber
	&= \Tr \bigg(  \left((I-QQ^\top) X - (I-QQ^\top)X(I-QQ^\top) \right) \\ \nonumber
	& \hspace{4cm} \times  \left(X(I-QQ^\top)  - (I-QQ^\top)X(I-QQ^\top) \right) \bigg)\\ \nonumber
	&\marksymb{=}{\ref{pfeq: facts 2}} \Tr((I-QQ^\top) X^2 (I-QQ^\top) - (I-QQ^\top) X(I-QQ^\top) X(I-QQ^\top) )\\ \nonumber
	&=\|X(I-QQ^\top)\|_2^2 - \| (I-QQ^\top )X (I-QQ^\top) \|_2^2\\ \label{pfeq: projection decreases norm}
	& \le \|X(I-QQ^\top)\|_2^2 = \|(X(I-QQ^\top))^\top \|_2^2 = \|(I-QQ^\top)X\|_2^2.
	\end{align}
	Now combining Eq.~\eqref{pfeq: square frob} with Ineq.~\eqref{pfeq: projection decreases norm} we get
	\begin{align} \label{pfeq: ineq frob 2 bound}
	\|X-\hat{X}\|_2^2=\|(I-QQ^\top)X\|_2^2 + \|QQ^\top X(I-QQ^\top)\|_2^2 \le 2 \|(I-QQ^\top)X\|_2^2.
	\end{align}
	Finally, taking the square root of the left and right hand sides of Ineq.~\eqref{pfeq: ineq frob 2 bound} we deduce Ineq.~\eqref{ineq: F norm bound} and complete the proof. }
\end{proof}


Using $X \approx {\hat{X}:=} Q Q^\top X Q Q^\top$ we can cheaply compute an approximate {PSD projection of $X \in \S_n$ based on the eigen-decomposition of the much smaller matrix $Q^\top X Q \in \mathbb S_k$, where $k\ll n$ . The core idea is that if \( X \) is well-approximated by \( \hat{X} = Q Q^\top X Q Q^\top \), then intuitively their positive semidefinite parts are also close, i.e., \( X_+ \approx \hat{X}_+ \). This intuition will be made precise in Corollary~\ref{cor: error bound on rand projection}, which provides explicit error bounds.}
	
	 By Lemma~\ref{lem: interchange PSD proj} (found in the appendix) it follows  $\hat{X}_+=(Q Q^\top X Q Q^\top)_+=Q (Q^\top X Q)_+ Q^\top$. Hence, to compute $\hat{X}_+$ we find the eigen-decomposition of the smaller matrix $Q^\top X Q=\hat{V} \hat{D} \hat{V}^\top \in \mathbb S_k$. Then by setting $\hat{U}=Q \hat{V}$ it follows $X_+ \approx \hat{X}_+ = \hat{U} \hat{D}_+ \hat{U}^\top \in \S_n^+.$ This method for approximate PSD projection is made explicit in Alg.~\ref{alg:approx proj}. This approximate method requires only \( \mathcal{O}(k n^2) \) operations, corresponding to the cost of the randomized subspace projection and eigen-decomposition of a \( k \times k \) matrix (see Page~249 of~\cite{halko2011finding}).



\begin{figure}
	\centering
	{
		\includegraphics[width=0.7\linewidth, trim = {1cm 0.4cm 0.8cm 0.4cm}, clip]{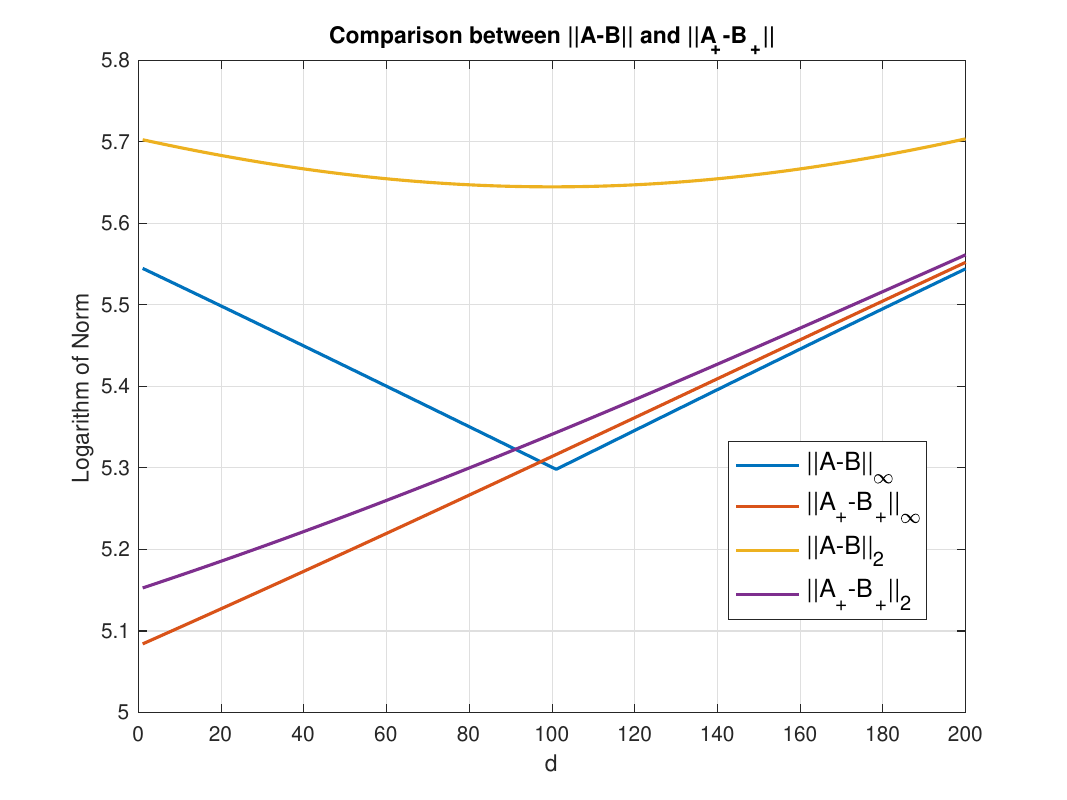}
		
	}
	\vspace{-10pt}
	\caption{\footnotesize Comparing $\|A-B\|$, $\|A_+-B_+\|$ where $A=\diag[1, 100]$ and $B = \begin{bmatrix}0 & 200\\ 200 & d \end{bmatrix}$ and $d$ is varied between $[0,200]$. Here, PSD projections are calculated using Eq.~\eqref{eq: analytical PSD proj using eigenvalue}, where the eigen-decomposition is calculated using MATLAB's \texttt{eig} function. It is clear that  {the Frobenius norm is non-expansive yet the bound,  $\|A_+-B_+\|_2 \le \|A-B\|_2$, may be loose for specific $A,B \in \S_n$. On the other hand, it} is not the case that spectral norm is non-expansive, i.e. $\|A_+-B_+\|_\infty \not \le \|A-B\|_\infty$ for all $d \in [0,200]$.} \label{fig: comparison norms} 
\end{figure}

\begin{algorithm}
	\caption{Randomized PSD Projection}\label{alg:approx proj}
	\textbf{Function} $\hat{X}_+=\texttt{ran$\_$proj}(X,k,l,q)$. \\
	\hspace*{\algorithmicindent} \textbf{Input:} $X\in \S_n$, {approximate rank} $k \in \N$, {oversampling parameter} $l \in \N$ \\
	\text{ } \hspace{1.7cm}  {and power iteration parameter $q \in \N$.} \\
	\hspace*{\algorithmicindent} \textbf{Output:} {Approximate PSD projection } $\hat{X}_+ \in \S_n^+$. 
	\begin{algorithmic}[1]
		\State $Q=\texttt{range$\_$finder} (X,k,l,q)$ \Comment{Alg.~\ref{alg:RSNA}}
		\State $[{\hat{V}},{\hat{D}}]=\texttt{eig}(Q^\top X Q)$ \Comment{$Q^\top X Q ={\hat{V}\hat{D} \hat{V}^T}$}
		\State $\hat{X}_+=Q{\hat{V} \hat{D}_+ \hat{V}^\top} Q^\top $ \Comment{{$\hat{X}:=QQ^\top X QQ^\top$}}
	\end{algorithmic}	
\end{algorithm}
\ifthenelse{\boolean{longver}}{}{
	\vspace{-0.2cm}}
We now provide an error bound on the quality of the approximate PSD projection resulting from Alg.~\ref{alg:approx proj}. To derive this error bound, we establish  Lipschitz-like bounds of a projection using the Frobenius and spectral norms. We will see that in the case of the Frobenius norm, the PSD projection bound presented later in Ineq.~\eqref{ineq: F non expansive projections}, is a matrix analogue of the well-established non-expansive property of projections onto closed convex subsets of $\mathbb{R}^n$, as depicted in Fig.~\ref{subfig: non expansive} (see Prop.~3.1.3 in~\cite{hiriart2004fundamentals}). This analogue result arises from the isomorphism between the space of symmetric matrices, $\S_n$, and the Euclidean space $\R^{n(n+1)/2}$ (illustrated for $\S_2$ in Fig.~\ref{fig: PSD proj}).  {However, this non-expansive property does not generalize to all matrix norms. In particular, it is known that the PSD projection is not non-expansive under the spectral norm~\cite[p.~184]{goulart2020accuracy}. We provide an explicit numerical counterexample illustrating this behaviour in Fig.~\ref{fig: comparison norms}, which also demonstrates how loose the Frobenius bound can be for specific matrices $A,B \in \S_n$. While these bounds may not be tight in all cases, they serve to illustrate how different problem configurations and parameters influence the overall projection error.}

Consequently, we instead provide two bounds for the difference in PSD projections under the spectral norm. The first bound, given in Ineq.~\eqref{ineq: spectral non expansive projections 1}, involves a constant that grows with the dimension of the matrices involved, whereas the second bound, given in Ineq.~\eqref{ineq: spectral non expansive projections 2}, provides an expansion constant involving the spectral norm magnitudes of the matrices involved.

\begin{prop}[PSD projection bounds] \label{prop: non expansive projections} 
	Consider symmetric matrices $X,Y \in \S_n$. Then the following inequalities hold,
	\ifthenelse{\boolean{longver}}{}{\vspace{-0.1cm}}
	\begin{align} \label{ineq: F non expansive projections}
		&	\|X_+-Y_+\|_2 \le \|X-Y\|_2.
		\\ \label{ineq: spectral non expansive projections 1}
		& \|X_+-Y_+\|_\infty \le \frac{1}{2} \left( \sqrt{n}+1 \right) \|X-Y\|_\infty.\\ \label{ineq: spectral non expansive projections 2}
		& 	\|X_+-Y_+\|_\infty \le  
		 \left( \frac{1}{2} +\frac{2}{\pi}  + \frac{1}{\pi} \log \left( \frac{\|X\|_\infty +\|Y\|_\infty}{\|X-Y\|_\infty}\right) \right)\|X-Y\|_\infty. 
	\end{align}
	

\end{prop} 
\ifthenelse{\boolean{longver}}{}{\vspace{-0.1cm}}
\begin{proof}
	We begin with an observation that given a symmetric matrix, $X\in \S_n$, we can construct the PSD projection, {$X_+$ given in Eq.~\eqref{opt: projecting onto PSD}}, using the polar decomposition of $X$ as follows,
	\ifthenelse{\boolean{longver}}{}{\vspace{-0.1cm}}\begin{align} \label{eq: projection in polar form}
		X_+=\frac{X+(X^\top X)^{\frac{1}{2}}}{2}.
	\end{align} 
 Eq.~\eqref{eq: projection in polar form} follows by considering eigen-decomposition of ${X}={U} {D} {U}^\top$ and then noting that $\frac{X+(X^\top X)^{\frac{1}{2}}}{2}=U \diag\left(\frac{D_{1,1} + |D_{1,1}|}{2}, \dots, \frac{D_{n,n} + |D_{n,n}|}{2}\right)U^\top={U} {D}_+ {U}^\top$, using the fact that $\max\{x,0\}=\frac{x+|x|}{2}$.  {Note that, while $X\in \S_n$ is a symmetric matrix (making $\sqrt{X^\top X}$ mathematically equivalent to $\sqrt{X^2}$), we use the redundant transpose notation to remind the reader that $\sqrt{X^\top X} \ne X$ unless $X \in \S_N^+$ is positive semi-definite.}
 \ifthenelse{\boolean{longver}}{Now it follows from Eq.~\eqref{eq: projection in polar form} and the triangle inequality that,
		\begin{align} \nonumber
			\|X_+  -Y_+\|& = \frac{1}{2} \|X-Y + (X^\top X)^{\frac{1}{2}} - (Y^\top Y)^{\frac{1}{2}}\|\\ \label{ineq: sub add matrix norm distance between PSD proj}
			& \le  \frac{1}{2} \left( \|X-Y\| + \|(X^\top X)^{\frac{1}{2}} - (Y^\top Y)^{\frac{1}{2}}\| \right),
		\end{align}
		where $\|\cdot\|$ is any valid matrix norm.
		
		To deduce Ineq.~\eqref{ineq: F non expansive projections} we simply apply Lemma~\ref{lem: An Inequality for the Frobenius norm} {(found in the Appendix)}, $\| \left(A^\top A \right)^{\frac{1}{2}}   - \left(B^\top B \right)^{\frac{1}{2}} \|_2 \le \| A -B\|_2$, to Ineq.~\eqref{ineq: sub add matrix norm distance between PSD proj} when the matrix norm is the Frobenious norm.}{
		\vspace{-0.2cm} 	To deduce Ineq.~\eqref{ineq: F non expansive projections} we simply apply Lemma 5.2 in~\cite{araki1971quasifree} to Ineq.~\eqref{ineq: sub add matrix norm distance between PSD proj} when the matrix norm is the Frobenious norm. }
 Inequality~\eqref{ineq: spectral non expansive projections 1} follows from the equivalence between the spectral norm and Frobenius norm: $\|A\|_\infty \le \|A\|_2 \le \sqrt{n} \|A\|_\infty$ for any $A \in \S_n$. Then by Lemma~\ref{lem: An Inequality for the Frobenius norm} and~\eqref{ineq: sub add matrix norm distance between PSD proj} we have that,
	\begin{align*}
		\|X_+-Y_+\|_\infty & \le \frac{1}{2} \left( \|X-Y\|_\infty + \|(X^\top X)^{\frac{1}{2}} - (Y^\top Y)^{\frac{1}{2}}\|_\infty \right)\\
		& \le \frac{1}{2} \|X-Y\|_\infty  + \frac{1}{2} \|(X^\top X)^{\frac{1}{2}} - (Y^\top Y)^{\frac{1}{2}}\|_2\\
		& \le \frac{1}{2} \|X-Y\|_\infty  + \frac{1}{2} \|X-Y\|_2 \le  \frac{1}{2} (1+ \sqrt{n})\|X-Y\|_\infty.
	\end{align*}	
	 Inequality~\eqref{ineq: spectral non expansive projections 2} follows by applying  Lemma~\ref{lem: An Inequality for the spectral norm} to Ineq.~\eqref{ineq: sub add matrix norm distance between PSD proj}.
\end{proof} 
\noindent If we assume $X$ and $Y$ are non-singular then we may improve the spectral norm bounds of Ineqs~\eqref{ineq: spectral non expansive projections 1} and~\eqref{ineq: spectral non expansive projections 2} by following the same argument of applying Ineq.~\eqref{ineq: sub add matrix norm distance between PSD proj} but now using {inequality VII.36 from~\cite{bhatia2013matrix}}\footnote{{This inequality also appears in~\cite{bhatia2010modulus}. Moreover, using the same proof strategy we may derive a bound for $\|X_+-Y_+\|$ with respect to the Schatten 1-norm from the Powers–Størmer inequality (Lem.~4.1 in~\cite{powers1970free}).}} for $\|(X^\top X)^{\frac{1}{2}} - (Y^\top Y)^{\frac{1}{2}}\|_\infty$.

Using the bounds from Proposition~\ref{prop: non expansive projections} we next bound the expected error of using Alg.~\ref{alg:approx proj} in both the Frobenius and Spectral norms. These error bounds will involve the following terms, defined for some $X \in \S_n$ and $k,l,q \in \N$,
\begin{align} \label{eps1}
	\eps_1(\{\sigma_i\}_{i=1}^n,k,l)&:=\sqrt{ \hspace{-4pt}\left( \hspace{-4pt}1+\frac{k}{l-1}  \hspace{-2pt}\right) \hspace{-4pt} \left(\sum_{j>k} \sigma_j^2  \right)}.\\ \label{eps2}
	\eps_2(\sigma,k,l,q)&:=	\hspace{-0.1cm}  \left( \hspace{-0.1cm} 1\hspace{-0.1cm} + \hspace{-0.1cm} \sqrt{\frac{k}{l-1}} \hspace{-0.1cm} + \hspace{-0.1cm}  \frac{e \sqrt{k+l}}{l}\sqrt{n-k} \right)^{\hspace{-0.15cm} \frac{1}{2q+1}} \hspace{-0.5cm} \sigma.
\end{align}


\ifthenelse{\boolean{longver}}{}{
	\vspace{-0.2cm}}
\begin{cor} \label{cor: error bound on rand projection}
	Consider $X\in \S_n$. Choose a target rank $k  \ge  2$ and an oversampling parameter $l \ge  2$, where $k + l \le n$. Execute Alg.~\ref{alg:approx proj} with $q=0$ to get an output of $\hat{X}_+ \in \S_n$, then
	{\begin{align} \label{ineq: F proj error}
			&	\mathbb{E}[	\|X_+ - \hat{X}_+\|_2] \le   {\sqrt{2}} \eps_1(\{\sigma_i(X)\}_{i=1}^n,k,l).
	\end{align} }
	\noindent Execute Alg.~\ref{alg:approx proj} for some $q\in \N$ to get an output of $\hat{X}_+ \in \S_n$, then
	\begin{align} \label{ineq: spectral proj error}  
			\mathbb{E}[	\|X_+ - \hat{X}_+\|_\infty] & \le \min \bigg\{(1+\sqrt{n}) \eps_2(\sigma_{k+1}(X),k,l,q),\\ \nonumber
		& \quad \left( 1 +\frac{4}{\pi}  + \frac{2}{\pi}(\log(2 \sigma_1(X)) )  \right)\eps_2(\sigma_{k+1}(X),k,l,q) \\ \nonumber
  & \hspace{2.5cm} +\frac{1}{\pi}\min\{e^{-1}, \sqrt{2\eps_2(\sigma_{k+1}(X),k,l,q)} \}\bigg\},
	\end{align}
	where the expectations of Ineqs~\eqref{ineq: F proj error} and~\eqref{ineq: spectral proj error} {are} taken with respect to the randomly generated matrix $\Omega$ and the $\eps_i$ terms can be found in Eqs~\eqref{eps1} and~\eqref{eps2}. 
\end{cor}
\ifthenelse{\boolean{longver}}{}{
	\vspace{-0.2cm}}
\begin{proof}
	{Execute Alg.~\ref{alg:approx proj} to get $Q \in \R^{n \times (k+l)}$ and output $\hat{X}_+$. Defining $\hat{X}:=QQ^\top X QQ^\top$ it follows by Lemma~\ref{lem: interchange PSD proj} that the PSD projection of $\hat{X}$ is equal to the output of Alg.~\ref{alg:approx proj}, that is $\hat{X}_+=(QQ^\top X QQ^\top)_+=Q(Q^\top X Q)_+Q^\top$. }
	
	 {By Ineqs~\eqref{ineq: F norm bound}} and~\eqref{ineq: F non expansive projections} we obtain
	\begin{align*}
		\|X_+ - \hat{X}_+\|_2 & \le \|X-\hat{X}\|_2  = \| X -Q Q^\top X Q Q^\top\|_2 \le  {\sqrt{2}} \|X-Q Q^\top X \|_2.
	\end{align*}
	\ifthenelse{\boolean{longver}}{\noindent Ineq.~\eqref{ineq: F proj error} then follows by applying Theorem~\ref{thm: expected error bounds of proto alg} in the appendix.}{
		Ineq.~\eqref{ineq: F proj error} the follows by applying Thm.~10.6 from~\cite{halko2011finding}.
	}
	
	We now show Ineq.~\eqref{ineq: spectral proj error}. From Ineqs~\eqref{ineq: 2 norm bound} and~\eqref{ineq: spectral non expansive projections 1} we have
	\begin{align} \nonumber
		\|X_+ - \hat{X}_+\|_\infty & \le\frac{1}{2} \left( \sqrt{n}+1 \right)\| X -Q Q^\top X Q Q^\top\|_\infty\\ \label{pfineq: inf bound 1}
		&\le \left( \sqrt{n}+1 \right)\| X -Q Q^\top X \|_\infty .
	\end{align}
	
	Moreover, from Ineq.~\eqref{ineq: spectral non expansive projections 2} we have that
	\begin{align} \nonumber 
		 \|X_+ - \hat{X}_+\|_\infty & \le 	\left( \frac{1}{2} +\frac{2}{\pi}  + \frac{1}{\pi} \log \left( \frac{\|X\|_\infty +\|\hat{X}\|_\infty}{\|X-\hat{X}\|_\infty}\right) \right)\|X-\hat{X}\|_\infty\\ \nonumber
		&= \left( \frac{1}{2} +\frac{2}{\pi}  + \frac{1}{\pi} \log \left( {\|X\|_\infty +\|\hat{X}\|_\infty}\right) \right)\|X-\hat{X}\|_\infty \\ \label{pfineq: inf}
		&\qquad \qquad  - \frac{1}{\pi}\log(\|X-\hat{X}\|_\infty)\|X-\hat{X}\|_\infty.
	\end{align}
	We now bound the various terms in Ineq.~\eqref{pfineq: inf}. Firstly, by the sub-multiplicative property of Schatten norms and Eq.~\eqref{eq: QQ 2 norm} we have that $\|\hat{X}\|_\infty=\|Q Q^\top X Q Q^\top\|_\infty  \le \|Q Q^\top\|_\infty^2 \|X\|_\infty= \|X\|_\infty$ and hence by the monotonicity property of the logarithm it follows that,
	\begin{align} \label{pfineq: inf 2}
		\log \left( {\|X\|_\infty +\|\hat{X}\|_\infty}\right) \le \log(2 \sigma_1(X) ).
	\end{align}
Now, applying the log inequality in Ineq.~\eqref{ineq: xlogx bound} (found in the Appendix), Ineq.~\eqref{ineq: 2 norm bound} and the monotonicity property of the square root, it follows that 
	\begin{align} \nonumber
			-\log(\|X-\hat{X}\|_\infty)\|X-\hat{X}\|_\infty & \le \min \left\{e^{-1}, \sqrt{ \|X-\hat{X}\|_\infty} \right\}\\ \label{pfineq: inf 3}
			& \le \min \left\{e^{-1}, \sqrt{ 2 \|X-QQ^\top X\|_\infty} \right\}.
		\end{align}
	Applying the bounds in Ineqs~\eqref{ineq: 2 norm bound}, \eqref{pfineq: inf 2} and~\eqref{pfineq: inf 3} along with Jensen's inequality, and noting $f(x)=\min\{e^{-1},\sqrt{2x} \}$ is a concave function,  to Ineq.~\eqref{pfineq: inf} we get that
	\begin{align} \nonumber
	&	 \mathbb{E}[	\|X_+ - \hat{X}_+\|_\infty  ] \le  \left( 1 +\frac{4}{\pi}  + \frac{2}{\pi} \log \left( 2 \sigma_1(X) \right) \right)\mathbb{E}[	\|X-QQ^\top X\|_\infty ] \\  \nonumber & \hspace{2.5cm} + \frac{1}{\pi} \mathbb{E}[	\min \left\{e^{-1}, \sqrt{ 2 \|X-QQ^\top X\|_\infty} \right\}] \\ \label{pfineq: inf bound 2}
		 & \le \left( 1 +\frac{4}{\pi}  + \frac{2}{\pi} \log \left( 2 \sigma_1(X) \right) \right)\mathbb{E}[	\|X-QQ^\top X\|_\infty ] \\ \nonumber 
		 & \hspace{2.5cm}+ \frac{1}{\pi} \min \left\{e^{-1}, \sqrt{ 2 \mathbb{E}[	\|X-QQ^\top X\|_\infty]} \right\}.
	\end{align}

	Applying the monotincity of the square root function and Ineq.~\eqref{ineq: proto error bound in spectral radius} (from Thm.~\ref{thm: expected error bounds of proto alg} in the appendix) to both Ineqs~\eqref{pfineq: inf bound 1} and~\eqref{pfineq: inf bound 2} and taking the minimum of these inequalities we deduce Ineq.~\eqref{ineq: spectral proj error}.
\end{proof}

{Interpreting the error bounds of Corollary~\ref{cor: error bound on rand projection}, the PSD projection accuracy of Alg.~\ref{alg:approx proj} is primarily controlled by the decay of singular values beyond rank $k$. The error scales with the corresponding optimal  {low rank matrix approximation} error from the Eckart-Young theorem \cite{eckart1936approximation} ($\sum_{j>k}\sigma_j^2(X)$ for the Frobenius norm or $\sigma_{k+1}(X)$ for the spectral norm), multiplied by factors that are improved by the choice of oversampling $l$ and power iterations $q$, yielding near-optimal approximations.}

\subsection{\textbf{A Scaled Randomized PSD Projection}} Real symmetric matrices, being a special case of \textit{normal matrices}, have the property that the absolute values of their eigenvalues are exactly equal to their singular values, i.e., $\{|\lambda_1(X)|,\dots,|\lambda_n(X)| \}=\{\sigma_1(X),\dots,\sigma_n(X)\}$. {As a result, if we first compute the optimal rank-$k$ approximation of a symmetric matrix and then project it onto the PSD cone, the outcome generally differs from the optimal rank-$k$ approximation of the PSD projection of the original matrix. The order of approximation and projection matters: computing a low-rank approximation before projection may retain singular values corresponding to negative eigenvalues, which are subsequently zeroed out during the PSD projection.} To see this consider the following example.

\begin{cex} \label{cex: optimal k rank matrix gives poor PSD proj} 
	Consider $X = \mathrm{diag}([-3,-2,1]).$
	{By the Eckart-Young theorem~\cite{eckart1936approximation}}, the optimal rank-1 approximation, in both the spectral and Frobenius norms, of $X$ is {obtained by retaining the component corresponding to the largest singular value. For this matrix, the largest singular value is $3$ (corresponding to the $-3$ entry),} yielding the optimal rank-1 approximation as $\hat{X}{:=}diag([-3,0,0])$. However $\hat{X}_+=0_{3\times 3}$. Alternatively the rank-1 matrix $\tilde{X}=diag([0,0,1])$ is a worse approximation of $X$ but has the same PSD projection as $X$: $\tilde{X}_+=X_+$.
\end{cex}

{Counterexample~\ref{cex: optimal k rank matrix gives poor PSD proj} illustrates a fundamental weakness of Alg.~\ref{alg:approx proj}: Given a matrix, $X \in \S_n$, it constructs a low rank approximate matrix, $X \approx \hat{X} =Q Q^\top X Q Q^\top$, and then projects this low rank matrix to get an approximate PSD projection, $\hat{X}_+$. When the dominant singular directions correspond to negative eigenvalues, the method preserves useless eigenvector directions in the initial low rank matrix approximation, which then leads to increased PSD projection error. To address this, we next introduce a scaling procedure that ensures the dominant singular values correspond to the positive eigenvalues.}

\ifthenelse{\boolean{longver}}{}{
	\vspace{-0.2cm}}
\begin{lem} \label{lem: singular values of B}
	Consider $X \in \S_n$. Suppose $\alpha:={|\min_i \lambda_i(X)|>0}$. Then \ifthenelse{\boolean{longver}}{$B:=\frac{X +\alpha I}{\alpha}$}{
		$B:=(X +\alpha I)/\alpha$} is a {PSD} symmetric matrix whose singular values and eigenvalues are such that
	\begin{align} \label{eq: eigen and singular values of B}
		\sigma_i(B)=\lambda_i(B)=\frac{\lambda_i(X) + \alpha}{\alpha}.
	\end{align}
	Moreover, {$\lambda_i(X)<0$ iff $\sigma_i(B) \in [0,1)$ and $\lambda_i(X) \ge 0$ iff $\sigma_i(B) \in [1,\infty)$.}
\end{lem} \ifthenelse{\boolean{longver}}{}{
	\vspace{-0.15cm}}
\ifthenelse{\boolean{longver}}{
	\begin{proof}
		{$X \in \S_n$ is a symmetric matrix and hence diagonalizable. Let $X=UDU^\top$ be the eigen-decomposition of $X$. It follows that $B:=(X +\alpha I)/\alpha$ is simultaneously diagonalizable with $B=U(D/\alpha+ I)U^\top$. Therefore the eigenvalues of $B$ are equal to $\frac{\lambda_i(X) + \alpha}{\alpha}$. Moreover, since $\alpha:=|\min_i \lambda_i(X)|$ it follows that all of the eigenvalues are positive and therefore $B \in \S_n^+$. The singular values and eigenvalues of PSD matricies are {identical}, thus proving Eq.~\eqref{eq: eigen and singular values of B}. Finally, {it follows immediately from Eq.~\eqref{eq: eigen and singular values of B} that $\lambda_i(X)<0$ iff $\sigma_i(B) \in [0,1)$ and $\lambda_i(X) \ge 0$ iff $\sigma_i(B) \in [1,\infty)$.}   }
	\end{proof}
}{
	\begin{proof}
		Uses a similar argument to Eq.~(6) in~\cite{francisco2017fixed}.
\end{proof}}

%

If we approximate \( B := \frac{X + \alpha I}{\alpha} \) by \( {B \approx \hat{B} :=Q Q^\top B Q Q^\top} \), we capture information about the leading singular values of \( B \), which correspond to the dominant positive eigenvalues of \( X \) (by Lemma~\ref{lem: singular values of B}). This enables us to retain the most important information for an approximate PSD projection.
Moreover, this scaling is particularly well suited to the power iteration method. As shown in Lemma~\ref{lem: singular values of B}, the negative eigenvalues of \( X \) correspond to the singular values of \( B \) in the interval \([0, 1]\), which are progressively compressed to zero as \( q \) increases in Alg.~\ref{alg:RSNA}. After approximating $B$ by $\hat{B}$, we can then project and re-scale $\hat{B}$ to derive an approximation of $X_+$. This method is summarized in Alg.~\ref{alg: scalled approx proj}. We now derive error bounds for computing an approximate PSD projection using Alg.~\ref{alg: scalled approx proj}. To do this, we first prove a key preliminary result that expresses the approximate PSD projection from Alg.~\ref{alg: scalled approx proj} in terms of \( B \).

\begin{algorithm}
	\caption{Scaled Randomized PSD Projection}\label{alg: scalled approx proj}
	\textbf{Function} $\hat{X}_+=\texttt{ran$\_$proj$\_$scal}(X,k,l,q,\alpha)$. \\
	\hspace*{\algorithmicindent} \textbf{Input:} $X\in \S_n$, {approximate rank} $k \in \N$, {oversampling parameter $l \in \N$ },\\
		\text{ } \hspace{1.7cm}  { power iteration parameter $q \in \N$.}
{and scaling parameter} $\alpha>0$. \\
	\hspace*{\algorithmicindent} \textbf{Output:} {Approximate PSD projection } $\hat{X}_+ \in \S_n^+$. 
	\begin{algorithmic}[1]
		\State $B=(X +\alpha I)/\alpha $
		\State $Q=\texttt{range$\_$finder} (B,k,l,q)$ \Comment{Alg.~\ref{alg:RSNA}}
		\State ${[\hat{V}},{\hat{D}}]=\texttt{eig}(Q^\top B Q)$ \Comment{$Q^\top B Q ={\hat{V} \hat{D} \hat{V}^T}$}
		\State $\hat{X}_+=\alpha Q {\hat{V}} (\max({\hat{D}},1) - I) {\hat{V}^\top} Q^\top $
	\end{algorithmic}	
\end{algorithm} 

\begin{lem} \label{lem: preliminary sclaed projection identity}
    Consider $X\in \S_n$, $k  \ge  2$, $k + l \le n$, $q \in \N$ and $\alpha: =|\min_i \lambda_i(X)|>0$. Execute the scaled randomized PSD algorithm (Alg.~\ref{alg: scalled approx proj}) to get an output of $\hat{X}_+ \in \S_n$ then
    \begin{align} \label{pfeq: X plus approx formula}
        \hat{X}_+ & =\alpha(QQ^\top(B-I)QQ^\top)_+\\ \label{eq: Simplification of scaled PSD}
        & = {\alpha(QQ^\top BQQ^\top-I)_+},
    \end{align}
    where matrices $B$ and $Q$ are found in Lines 1 and 2 respectively in Alg.~\ref{alg: scalled approx proj}.
\end{lem}
\begin{proof}
 {We first show Eq.~\eqref{pfeq: X plus approx formula}}.    From Alg.~\ref{alg: scalled approx proj} we have $\hat{X}_+=\alpha Q\hat{V} (\max(\hat{D},1) - I) \hat{V}^\top Q^\top $ where $\hat{V}\hat{D}\hat{V}^\top = Q^\top BQ $ and $\hat{V} \hat{V}^\top=\hat{V}^\top\hat{V}= I$. Therefore, using $Q^\top Q=I$ and Lemma~\ref{lem: interchange PSD proj} (found in the Appendix) we have that
	\begin{align} \nonumber 
		\hat{X}_+ & = \alpha Q\hat{V} (\max(\hat{D},1) - I) \hat{V}^\top Q^\top  =\alpha Q\hat{V} (\hat{D} - I)_+ \hat{V}^\top Q^\top \\  \nonumber
		& \marksymb{=}{\ref{eq: interchange PSD projection}} \alpha Q(\hat{V} (\hat{D} - I) \hat{V}^\top)_+ Q^\top  = \alpha Q(Q^\top BQ-I)_+ Q^\top\\ \nonumber
		&=\alpha Q(Q^\top BQ-Q^\top Q)_+ Q^\top =\alpha Q(Q^\top (B-I)Q)_+ Q^\top\\ \nonumber
  &\marksymb{=}{\ref{eq: interchange PSD projection}} \alpha(QQ^\top(B-I)QQ^\top)_+.
	\end{align} 
	
 {We next show Eq.~\eqref{eq: Simplification of scaled PSD}. Let us denote $X_1:=\alpha QQ^\top(B-I)QQ^\top$ and $X_2:=\alpha(QQ^\top B QQ^\top -I)$. Showing Eq.~\eqref{eq: Simplification of scaled PSD} is equivalent to showing $(X_1)_+=(X_2)_+$. Note, by definition we have 
\begin{align} \label{pfeq: X1 and X2 relationship}
	X_2=X_1-\alpha (I - QQ^\top ).
\end{align}
Since $QQ^\top(I-QQ^\top)=QQ^\top-QQ^\top QQ^\top=QQ^\top-QQ^\top=0$ we also have
\begin{align} \label{pfeq: X1 null}
	(I-QQ^\top)X_1=\alpha(I-QQ^\top)QQ^\top(B-I)QQ^\top=0.
\end{align}
\begin{flalign}
	\makebox[0pt][l]{Similarly,} && X_1(I-QQ^\top)=0. && \label{pfeq: X1 null 2}
\end{flalign}
Eqs~\eqref{pfeq: X1 null} and \eqref{pfeq: X1 null 2} show that the column vectors of $(I-QQ^\top)$ are a subset of the null space of  $X_1$ and hence are also a subset of the null space of $X_1^2$. Since the null space of $X_1^2$ and $\sqrt{X_1^2}$ are equal (see Theorem 7.2.7 in~\cite{horn2012matrix}) it follows that,
\begin{align} \label{pfeq: sqrt X1 null 2}
	\sqrt{X_1^2}	(I-QQ^\top)=	(I-QQ^\top) \sqrt{X_1^2}	=0.
\end{align}

It then follows that 
\begin{align} \label{pfeq: X1 sqrt}
\sqrt{X_1^2 + \alpha^2 (I - QQ^\top )}=\sqrt{X_1^2} + \alpha(I - QQ^\top ),
\end{align}
because the square root of a PSD matrix, $A \in S_n^+$, is uniquely defined as the PSD matrix $B \in S_n^+$ such that $B^\top B=B^2=A$ is satisfied\footnote{We denote the square root of $A \in S_n^+$ by $\sqrt{A}:=B$ (see notational section).}. Thus, since $X_1^2$ is PSD, $\alpha^2 (I - QQ^\top )$ is PSD it follows $X_1^2 + \alpha^2 (I - QQ^\top ) \in S_n^+$ is PSD and hence the left hand side of Eq.~\eqref{pfeq: X1 sqrt} is uniquely defined. The right hand side of Eq.~\eqref{pfeq: X1 sqrt} is also PSD, so to show Eq.~\eqref{pfeq: X1 sqrt} it follows that all we need to do is to square the right hand side and shows that this equals $X_1^2 + \alpha^2 (I - QQ^\top )$:
\begin{align*}
& \left( \sqrt{X_1^2} + \alpha(I - QQ^\top ) \right)^2\\
& = \sqrt{X_1^2}\sqrt{X_1^2} + \alpha \sqrt{X_1^2}(I-QQ^\top) + \alpha (I-QQ^\top) \sqrt{X_1^2} +  \alpha^2 (I - QQ^\top )^2\\
& \marksymb{=}{\eqref{pfeq: facts 2} , \eqref{pfeq: sqrt X1 null 2}} X_1^2 + \alpha^2 (I - QQ^\top ).
\end{align*}

Now, using Eq.~\eqref{eq: projection in polar form} we have that
\begin{align*}
&	(X_2)_+ \marksymb{=}{\ref{eq: projection in polar form} }  \frac{1}{2} \left(X_2 + \sqrt{X_2^\top X_2} \right)\\
	& \marksymb{=}{\ref{pfeq: X1 and X2 relationship}} \frac{1}{2} \left( X_1-\alpha (I - QQ^\top ) + \sqrt{(X_1-\alpha (I - QQ^\top ))^\top (X_1-\alpha (I - QQ^\top ))} \right)\\
		& \marksymb{=}{\eqref{pfeq: facts 2} , \eqref{pfeq: X1 null}} \frac{1}{2} \left( X_1-\alpha (I - QQ^\top ) + \sqrt{X_1^2 + \alpha^2 (I - QQ^\top )} \right)\\
		&	\marksymb{=}{\ref{pfeq: X1 sqrt}} \frac{1}{2} \left( X_1-\alpha (I - QQ^\top ) + \sqrt{X_1^2} + \alpha (I - QQ^\top )\right) 	\marksymb{=}{\ref{eq: projection in polar form} } (X_1)_+,
\end{align*}
Hence, Eq.~\eqref{eq: Simplification of scaled PSD} is shown and the proof is complete. }
\end{proof}
We now use Lemma~\ref{lem: preliminary sclaed projection identity} to derive error bounds for the quality of the PSD projection given by Alg.~\ref{alg: scalled approx proj} in both the Frobenius and spectral norms.
\begin{prop} \label{prop: proj error 2}
	Consider $X\in \S_n$ and ${\alpha: =|\min_i \lambda_i(X)|>0}$. Choose a target rank $k  \ge  2$ and an oversampling parameter $l \ge  2$, where $k + l \le n$. Execute the scaled randomized PSD algorithm (Alg.~\ref{alg: scalled approx proj}) with $q=0$ to get an output of $\hat{X}_+ \in \S_n$. Then
	\ifthenelse{\boolean{longver}}{}{
		\vspace{-0.2cm}}	{\begin{align} \label{ineq: proj error 2}
		&	\mathbb{E}[	\|X_+  -  \hat{X}_+\|_2] \le  {\sqrt{2}} \eps_1(\{(\lambda_i(X) + \alpha)\}_{i=1}^n,k,l).
   \end{align}
   	\noindent Execute the scaled randomized PSD algorithm (Alg.~\ref{alg: scalled approx proj}) for some $q\in \N$ to get an output of $\hat{X}_+ \in \S_n$, then
    \begin{align} \nonumber
			&	\mathbb{E}[	\|X_+  -  \hat{X}_+\|_\infty] \le \min \bigg\{ (1+\sqrt{n}) \eps_2((\lambda_{k+1}(X)+\alpha),k,l,q) , \\ \nonumber
			& \qquad \left( 1 +\frac{4}{\pi}  + \frac{2}{\pi}(\log(2 \sigma_1(X)) )  \right)\eps_2((\lambda_{k+1}(X)+\alpha),k,l,q)\\ \label{ineq: proj error 3}
   & \qquad +\frac{1}{\pi}\min\{e^{-1},\sqrt{2 \eps_2((\lambda_{k+1}(X)+\alpha),k,l,q) } \} \bigg\}
	\end{align} }
	where the expectation is taken with respect to the randomly generated matrix $\Omega$ and the $\eps_i$ terms can be found in Eqs~\eqref{eps1} and~\eqref{eps2}. 
\end{prop}

\begin{proof}
	Let $B:=(X +\alpha I)/\alpha$. Then since $\alpha>0$ and hence it follows that
	\begin{align} \label{X_+ in terms of B}
		X_+=\alpha(B-I)_+.
	\end{align}
	
	We split the proof into two parts. In part one we show our Frobenius norm error bound, given in Ineq.~\eqref{ineq: proj error 2}. In part two we show our spectral norm error bound given in Ineq.~\eqref{ineq: proj error 3}.

\textbf{Part One:}	
	\begin{align} \label{ineq: PSD proj prelim bound}
			\|X_+-\hat{X}_+\|_2 & \marksymb{=}{\ref{X_+ in terms of B}} \|\alpha(B-I)_+-\hat{X}_+\|_2 \\ \nonumber
   & \marksymb{=}{\ref{eq: Simplification of scaled PSD}} \| \alpha (B-I)_+ - \alpha  {(QQ^\top B QQ^\top-I)_+  } \|_2   \\ \nonumber
   &\le \alpha \| (B-I)_+ -  {(QQ^\top B QQ^\top-I)_+  }   \|_2\\ \nonumber
		&  \marksymb{\le}{\ref{ineq: F non expansive projections}} \alpha \| (B-I) -  {(QQ^\top B QQ^\top-I) }   \|_2 \\ \nonumber
		&=\alpha \| (B - (QQ^\top B QQ^\top))   \|_2  \\ \nonumber
		& \marksymb{\le}{\ref{ineq: F norm bound}}  \alpha  {\sqrt{2}} \| B - QQ^\top B \|_2 
	\end{align}
	We now take the expectation  {of Ineq.~\eqref{ineq: PSD proj prelim bound} and bound} $\alpha \mathbb{E}[\|(I-QQ^\top)B\|_2]$. 

	\begin{align}  \nonumber
		&  \alpha \mathbb{E}[\|(I-QQ^\top)B\|_2]  \marksymb{\le}{\ref{ineq: randomized frobenious error bound}} \alpha \sqrt{\left(1+\frac{k}{l-1} \right) \left(\sum_{j>k} \sigma_j(B)^2  \right)}\\ 
		& \marksymb{=}{\ref{eq: eigen and singular values of B}} \alpha \sqrt{\left(1+\frac{k}{l-1} \right) \left(\sum_{j>k} \left(\frac{\lambda_j(X) + \alpha}{\alpha}\right)^2  \right)}  
 \label{ineq: I QQ B bound}
		\marksymb{=}{\ref{eps1}} \eps_1(\{(\lambda_i(X) + \alpha)\}_{i=1}^n,k,l).
	\end{align}
	Combining bounds from Ineqs~\eqref{ineq: PSD proj prelim bound} and~\eqref{ineq: I QQ B bound} we get Ineq.~\eqref{ineq: proj error 2}.
	
\textbf{Part Two:}	We now show Ineq.~\eqref{ineq: proj error 3}. First note by Eq.~\eqref{eq: QQ eigenvalues} that the eigenvalues of $QQ^\top$ are either zero or one, implying $\| I -QQ^\top   \|_\infty \le 1$. Now, using the triangle inequality along with Ineqs~\eqref{ineq: 2 norm bound} and~\eqref{ineq: spectral non expansive projections 1} we get that,
	\begin{align} \label{ineq: PSD proj prelim bound 2}
			\mathbb{E} \left[\|X_+-\hat{X}_+\|_\infty  \right] & \marksymb{=}{\ref{eq: Simplification of scaled PSD}}  { \mathbb{E} \left[\| \alpha(B-I)_+ - \alpha (QQ^\top B QQ^\top-I)_+   \|_\infty \right] } \\ \nonumber  & \le \alpha \mathbb{E} \left[\| (B-I)_+ -  {(QQ^\top B QQ^\top -I)_+  } \|_\infty \right]\\ \nonumber
		& \marksymb{\le}{\ref{ineq: spectral non expansive projections 1}} \frac{\alpha}{2}(1+\sqrt{n}) \mathbb{E} \left[\| (B-I) -  {(QQ^\top B QQ^\top -I)  }   \|_\infty \right] \\ \nonumber
		&  {=} { \frac{\alpha}{2}(1+\sqrt{n}) \mathbb{E} \left[ \| (B - (QQ^\top B QQ^\top))\|_\infty \right] }\\ \nonumber
		& \marksymb{\le}{\ref{ineq: 2 norm bound}} \alpha (1+\sqrt{n})  \mathbb{E} \left[\| B - QQ^\top B \|_\infty \right] .
	\end{align}

By Eq.~\eqref{eq: Simplification of scaled PSD} we have that $\hat{X}_+= \alpha( {QQ^\top B QQ^\top-I})_+$ and hence $\hat{X}= \alpha(QQ^\top B QQ^\top-I)=QQ^\top (X +\alpha I) QQ^\top -\alpha I=QQ^\top X QQ^\top +\alpha(QQ^\top-I)$. Hence, using Lemma~\ref{lem: spectal projection bound} (found in the appendix) and $(QQ^\top-I)^2=(QQ^\top-I)$ we get that,
 {\begin{align} \nonumber
	\|\hat{X}\|_\infty & = \| QQ^\top X QQ^\top +(QQ^\top-I)^\top \alpha I (QQ^\top-I)\|_\infty \\ \label{pfeq:bound hat X}
	& \marksymb{\le}{\ref{eq: inf norm bound}} \max \{ \|X\|_\infty, \alpha \} = \max \left\{ \sigma_1(X), \left|\min_i \lambda_i(X) \right| \right\} =\sigma_1(X). 
\end{align} }

Now, by a similar argument to Ineq.~\eqref{pfineq: inf bound 2} we get,
	\begin{align} \nonumber
		&	\mathbb{E}[\|X_+-\hat{X}_+\|_\infty] \marksymb{\le}{\eqref{ineq: spectral non expansive projections 2}, \eqref{eq: Simplification of scaled PSD} }	\mathbb{E} \bigg[ \alpha  \left( \frac{1}{2} +\frac{2}{\pi}  + \frac{1}{\pi} \log \left( \|X\|_\infty +\|  {\hat{X}} \|_\infty  \right) \right)\\ \nonumber
		& \hspace{1.5cm}  \times \|B-I-(QQ^\top BQQ^\top-I)\|_\infty  -\frac{1}{\pi}\log(\|X- \hat{X}\|_\infty)\|X- \hat{X}\|_\infty \bigg] \\ \nonumber
		& \marksymb{\le}{\ref{pfeq:bound hat X}} \mathbb{E} \bigg[ \alpha  \left( \frac{1}{2} +\frac{2}{\pi}  + \frac{1}{\pi} \log \left(2\sigma_1(X)  \right) \right) \| B - QQ^\top B Q Q^\top \|_\infty  \\ \nonumber
  & \qquad + \frac{1}{\pi} \min \left\{e^{-1},\sqrt{ \alpha \| B - QQ^\top B Q Q^\top \|_\infty  } \right\} \bigg] \\ \nonumber
		& \marksymb{\le}{\ref{ineq: 2 norm bound}}   2 \alpha  \left( \frac{1}{2} +\frac{2}{\pi}  + \frac{1}{\pi} \log \left(2\sigma_1(X)  \right) \right) \mathbb{E}[ \|(I - QQ^\top) B \|_\infty] \\ \label{ineq: PSD proj prelim bound 3}
  & \qquad+ \frac{1}{\pi}\min \left\{e^{-1},\sqrt{ 2\alpha \mathbb{E}[\| (I - QQ^\top) B  \|_\infty] } \right\}.
	\end{align}
	We next bound $\alpha \mathbb{E}[\|(I-QQ^\top)B\|_\infty]$ in Ineqs~\eqref{ineq: PSD proj prelim bound 2} and~\eqref{ineq: PSD proj prelim bound 3}:
 \begin{align} \label{ineq: I QQ B bound spectral}
      \alpha \mathbb{E}[ &\|(I-QQ^\top)B\|_\infty] \\ \nonumber
     & \marksymb{\le}{\ref{ineq: proto error bound in spectral radius}} \hspace{-0.1cm}  \left(1+\sqrt{\frac{k}{l-1}} + \frac{e \sqrt{k+l}}{l}\sqrt{\min\{m,n\}-k} \right)^{\frac{1}{2q+1}}  \alpha \sigma_{k+1}(B)\\ \nonumber
     &\marksymb{=}{\ref{eq: eigen and singular values of B}} \hspace{-0.1cm}  \left(1+\sqrt{\frac{k}{l-1}} + \frac{e \sqrt{k+l}}{l}\sqrt{\min\{m,n\}-k} \right)^{\frac{1}{2q+1}}  (\lambda_{k+1}(X) + \alpha)\\ \nonumber
     & \marksymb{=}{\ref{eps2}} \eps_2((\lambda_{k+1}(X)+\alpha),k,l,q).
 \end{align}
 {Where Ineq.~\eqref{ineq: proto error bound in spectral radius} is from~\cite{halko2011finding} and provided in the appendix}.
 Applying this to the bounds given in Ineqs~\eqref{ineq: PSD proj prelim bound 2} and~\eqref{ineq: PSD proj prelim bound 3} and taking the minimum of these bounds we derive Ineq.~\eqref{ineq: proj error 3}.
\end{proof}


	{In contrast to the bounds in Corollary~\ref{cor: error bound on rand projection}, the error bounds in Proposition~\ref{prop: proj error 2} are governed by the residual of the \textit{shifted} eigenvalues, $\lambda_i(X) + \alpha$, rather than the singular values, $\sigma_i(X)$. This shifted spectrum may exhibit faster decay.  {Consequently, unlike Corollary~\ref{cor: error bound on rand projection}, the error bounds of Proposition~\ref{prop: proj error 2} can be tighter than the optimal error predicted by the Eckart-Young theorem for low rank matrix approximation; consider the case $\lambda_i(X) + \alpha=0$ for all $i \ge k$ then Ineqs~\eqref{ineq: proj error 2} and~\eqref{ineq: proj error 3} imply zero expected PSD projection error.} 



  In the following corollary we give an example of a class of matrices where the right hand side of the bound given in Ineq.~\eqref{ineq: F proj error} is greater than right hand side of the bound given in Ineq.~\eqref{ineq: proj error 2}. This indicates that, for at least this class of matrix, it is more accurate to use the scaled randomized PSD projection (Alg.~\ref{alg: scalled approx proj}) than the randomized PSD projection without scaling (Alg.~\ref{alg:approx proj}). Note that we can also construct a similar example for the spectral norm but omit this for the sake of brevity. In Figure~\ref{fig: eig test} we plot the eigen and singular values for a particular instance of a matrix within this class for $n=1000$. {This figure visualizes how the leading singular values and eigenvalues may differ in order, highlighting a case of matrix where the second cluster of singular values correspond to a cluster of negative eigenvalue and thus provides no useful information for PSD projection.}

\begin{cor}
\label{cor: when sclaing bound is smaller}
		Consider the following matrix,
	\begin{align} \label{eq: matrix where scale is better}
X=Y^\top\mathrm{diag}([-\beta_1 \mathbf{1}_{n/4} 
		,-\beta_2 \mathbf{1}_{4/2},\beta_3 \mathbf{1}_{n/4} 
		,\beta_4 \mathbf{1}_{4/2}])Y \in \S_n,
	\end{align}
		where $Y \in \S_n$ is an orthogonal matrix and $\beta_3>\beta_1>\beta_4>\beta_2>0$ is such that $(\beta_1-\beta_2)^2 < \beta_2^2 + \beta_4^2$. Then
		\begin{align*}
		\eps_1(\{(\lambda_i(X)+\alpha)\}_{i=1}^n,n/2,l) <  \eps_1(\{\sigma_i(X)\}_{i=1}^n,n/2,l) ,
		\end{align*}
  with the $\eps_i$ terms from Eqs~\eqref{eps1} and~\eqref{eps2}.  {This indicates the right hand side of the bound for the scaled PSD projection error, given in  Ineq.~\eqref{ineq: proj error 2}, is less than the right hand side of the bound for the vanilla PSD projection error given in Ineq.~\eqref{ineq: F proj error} for $k=n/2$. }
\end{cor} 
\begin{proof}
	It follows that 
	\begin{align*}
		\sigma_i(X)=\begin{cases} \beta_3 \text{ } 1\le i \le n/4 \\ \beta_1 \text{ } n/4 < i \le n/2 \\
  \beta_4 \text{ } n/2 < i \le 3n/4 \\ 
   \beta_2 \text{ } 3n/4 < i \le n\end{cases} \lambda_i(X)=\begin{cases} \phantom{-} \beta_3 \text{ }  1 \le i \le n/4 \\ \phantom{-} \beta_4 \text{ } n/4< i \le n/2\\  -\beta_2 \text{ } n/2< i \le 3n/4\\
   -\beta_1 \text{ } 3n/4< i \le n\end{cases}  
	\end{align*}
	Hence $\alpha=\beta_1$ (the magnitude of the most negative eigenvalue). Now, 
	\begin{align*}
	\eps_1(\{(\lambda_i(X)+\alpha)\}_{i=1}^n,n/2,l) &= \sqrt{1+ \frac{n}{2(l-1)}}\sqrt{\frac{n}{4}(\beta_1-\beta_2)^2}\\
	& < \sqrt{1+ \frac{n}{2(l-1)}}\sqrt{\frac{n}{4}(\beta_2^2 + \beta_4^2)}\\
	& = 	\eps_1(\{(\sigma_i(X))\}_{i=1}^n,n/2,l) 
	\end{align*}
	\end{proof}


To execute Alg.~\ref{alg: scalled approx proj} we estimate $\alpha \approx |\min_i \lambda_i(X)|$. Fortunately, {one possible way\footnote{{Other possible methods exist for estimating extremal eigenvalues such as LOBPCG~\cite{knyazev2001toward}, Jacobi–Davidson Iteration~\cite{sleijpen2000jacobi}, sketching~\cite{swartworth2023optimal}, etc.}}  to calculate} such an estimation {is to use the cheap to compute}\footnote{{Each iteration is dominated by matrix vector multiplication with a well known FLOP count of $\mathcal{O}(n^2)$~\cite{hunger2005floating}.}} classical power iteration method that dates back to the 1920's~\cite{mises1929praktische,bai2021power}, see Alg.~\ref{alg: max singular} for the basic procedure\footnote{{See \cite{yang2020analysis} for a numerical stability analysis of this method}.}. The maximum singular value of a symmetric matrix, estimated using Alg.~\ref{alg: max singular}, either corresponds to the absolute value of the minimal or maximal eigenvalue. We are only interested in computing the absolute value of the minimal eigenvalue. Therefore, to find the minimal eigenvalue we must then make the following transformation: $Y=X-\sigma_{1}(X)I$. It follows that $|\min_i \lambda_i(X)|=|\sigma_{1}(Y)-\sigma_{1}(X)|$, as shown in Lemma~\ref{lem: caluclate min eigenvalue}. The entire procedure for calculating the absolute value of the minimal eigenvalue is summarized in Alg.~\ref{alg: min eig}.  
\begin{center}
	\begin{algorithm}
		\caption{Approximate maximum singular value}\label{alg: max singular}
		\textbf{Function} $\sigma_{1}=\texttt{power$\_$iteration}(X,N)$. \\
		\hspace*{\algorithmicindent} \textbf{Input:} $X\in \S_n$ {and number of power iterations} $N \in \N$. \\
		\hspace*{\algorithmicindent} \textbf{Output:} {Approximate largest singular value} $\sigma_{1}>0$. 
		\begin{algorithmic}[1]
			\State $v=\texttt{randn}(n,1) $ \Comment{Random initialization {with $v_i\marksymb{\sim}{iid}\mathcal N (0,1)$}}
			\For{$i=1:N$}
			\State $v=Xv/\| X v\|_2$ 
			\EndFor
			\State $\sigma_{1}=\| X v\|_2$
		\end{algorithmic}	
	\end{algorithm} 
\end{center}
\begin{center}
	\begin{algorithm}
		\caption{Approximate absolute value of minimum eigenvalue}\label{alg: min eig}
		\textbf{Function} $\lambda_{min}=\texttt{min$\_$eig}(X,N)$. \\
		\hspace*{\algorithmicindent} \textbf{Input:} $X\in \S_n$ {and number of power iterations} $N \in \N$. \\
		\hspace*{\algorithmicindent} \textbf{Output:} {Approximate absolute value of minimal eigenvalue} $\lambda_{min}>0$. 
		\begin{algorithmic}[1]
			\State $\sigma_1=\texttt{power$\_$iteration}(X,N) $  \Comment{Alg.~\ref{alg: max singular}}
			\State $Y=X-\sigma_1 I$
			\State $\sigma_2=\texttt{power$\_$iteration}(Y,N) $ \Comment{Alg.~\ref{alg: max singular}}
			\State $\lambda_{min}=|\sigma_1-\sigma_2|$
		\end{algorithmic}	
\end{algorithm} \end{center}

\begin{figure}
\subfloat[ {\scriptsize Eigenvalue clusters $\lambda_1(X)>\lambda_2(X)>0>\lambda_3(X)>\lambda_4(X)$.  } \label{subfig: eig plot}]{\includegraphics[width=0.49 \linewidth, trim ={1cm 0cm 1.25cm 0cm}, clip]{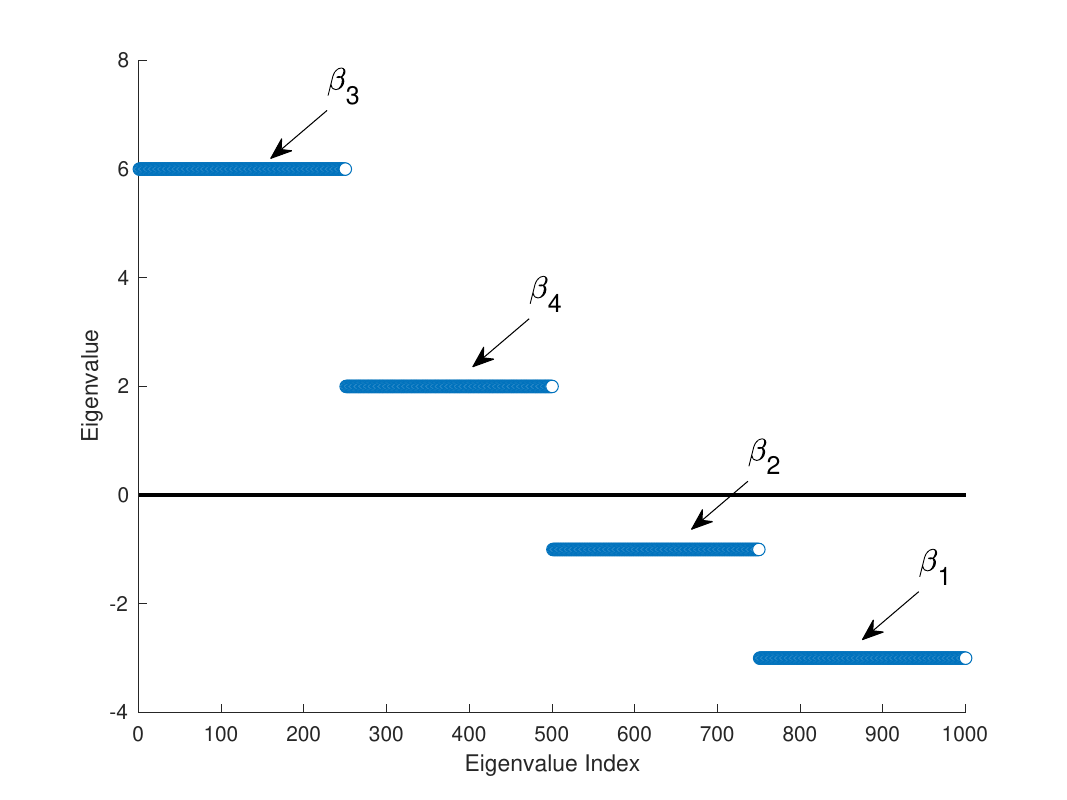}}\hfill 
		\subfloat[ {\scriptsize Singular value clusters $\sigma_1(X)>\sigma_2(X)>\sigma_3(X)>\sigma_4(X)>0$.} \label{subfig: sing test plot}]{\includegraphics[width=0.49 \linewidth, trim = {1cm 0cm 1.25cm 0cm}, clip]{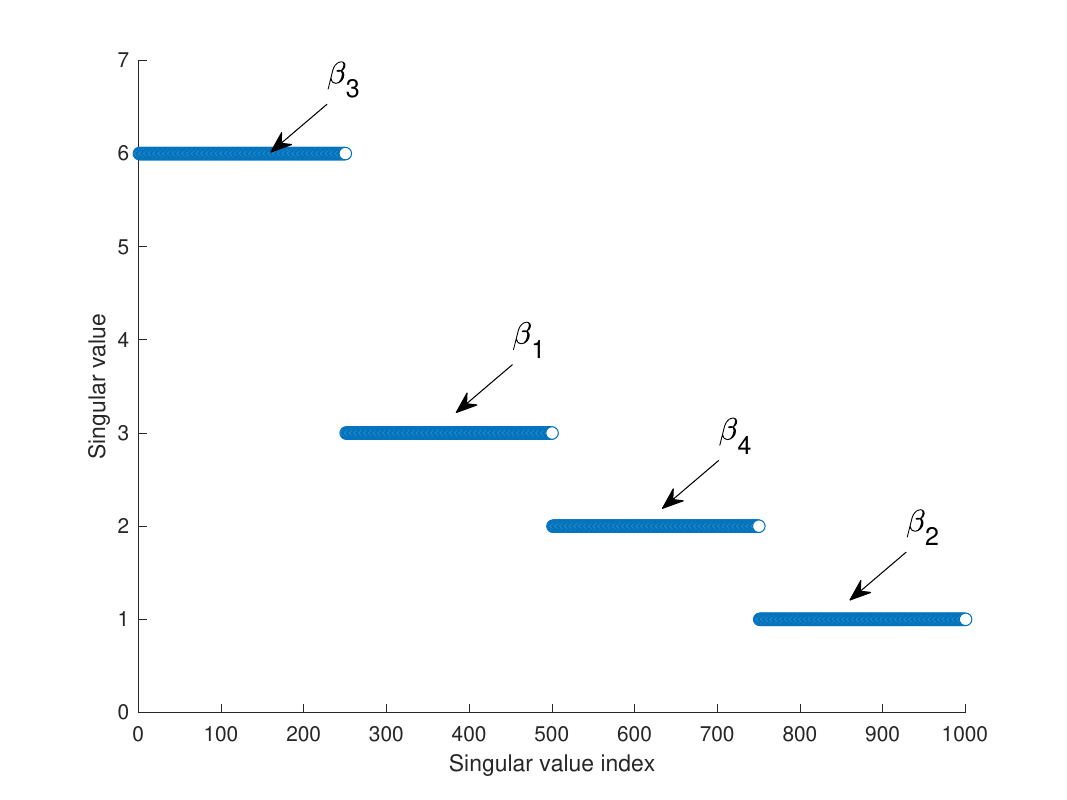}}\hfill 
		\caption{ \footnotesize Plot showing the distribution of eigenvalues and singular values of the matrix given in Eq.~\eqref{eq: matrix where scale is better} and corresponding associated scalars $\beta_3>\beta_1>\beta_4>\beta_2>0$.} \label{fig: eig test}
\end{figure}

\begin{figure}
		\subfloat[ {\scriptsize  Logarithm of numerical Frobenius error of Algs~\ref{alg:approx proj} and~\ref{alg: scalled approx proj}. \vspace{-10pt}} \label{subfig: actual_frob_error} ]{\includegraphics[width=0.49 \linewidth, trim = {1cm 0cm 1cm 0cm}, clip]{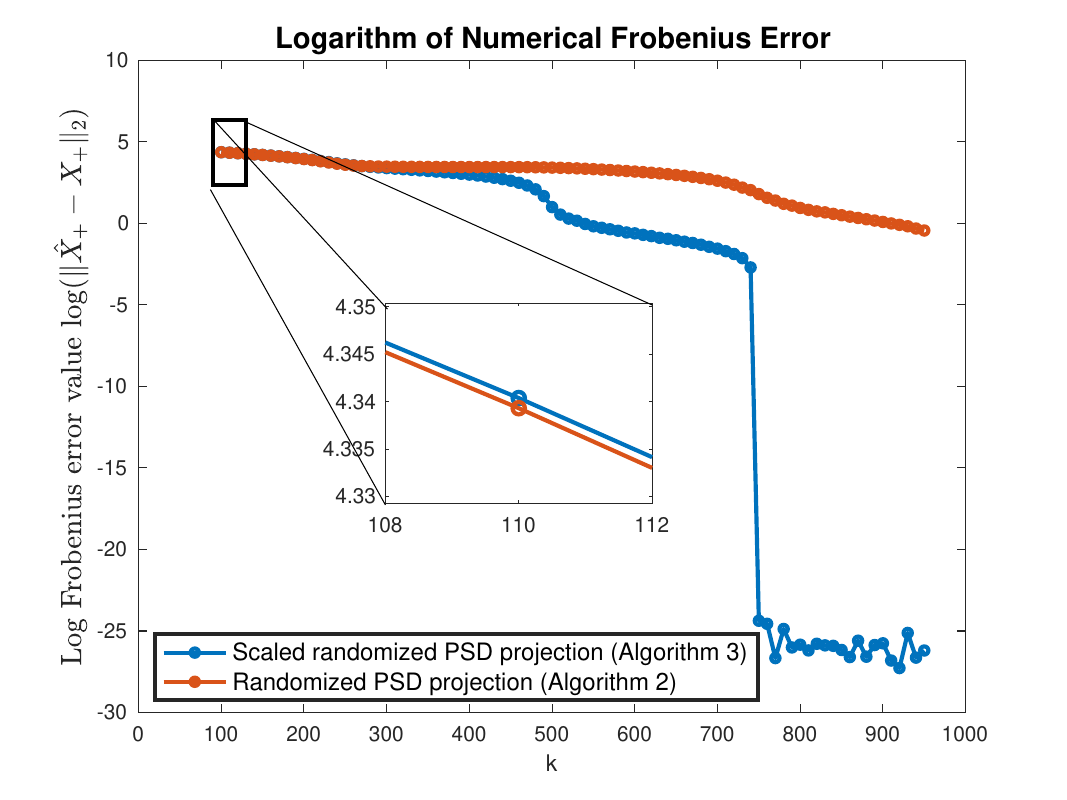}} \hfill 
		\subfloat[ {\scriptsize Logarithm of Frobenius error bounds of Ineqs~\eqref{ineq: F proj error} and~\eqref{ineq: proj error 2}. \vspace{-10pt}} \label{subfig: RHS_bound_frob_error}]{\includegraphics[width=0.49 \linewidth, trim = {0.75cm 0cm 1cm 0cm}, clip]{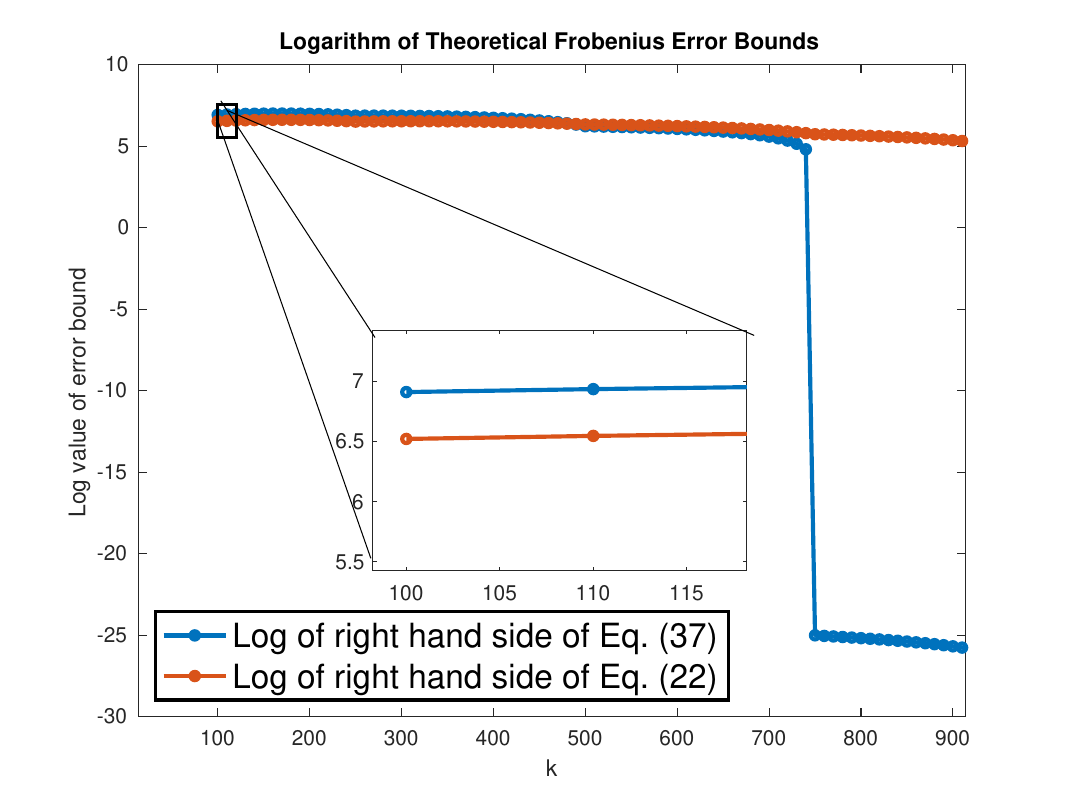}}\hfill 
		\subfloat[ {\scriptsize Numerical spectral error of Algorithms~\ref{alg:approx proj} and~\ref{alg: scalled approx proj}.  } \label{subfig: actual_spec_error}]{\includegraphics[width=0.49 \linewidth, trim ={1cm 0cm 1cm 0cm}, clip]{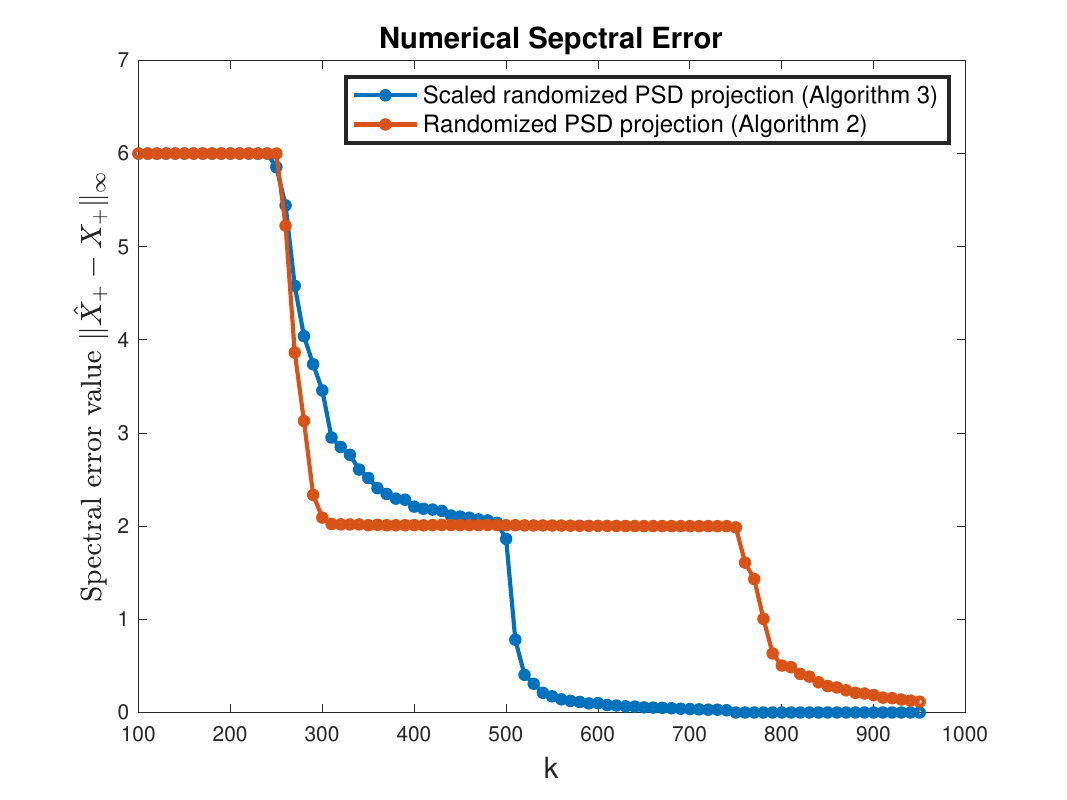}}\hfill 
		\subfloat[ {\scriptsize Spectral error bounds of Ineqs~\eqref{ineq: spectral proj error} and~\eqref{ineq: proj error 3}.} \label{subfig: RHS_bound_spec_error}]{\includegraphics[width=0.49 \linewidth, trim = {0.75cm 0cm 1cm 0cm}, clip]{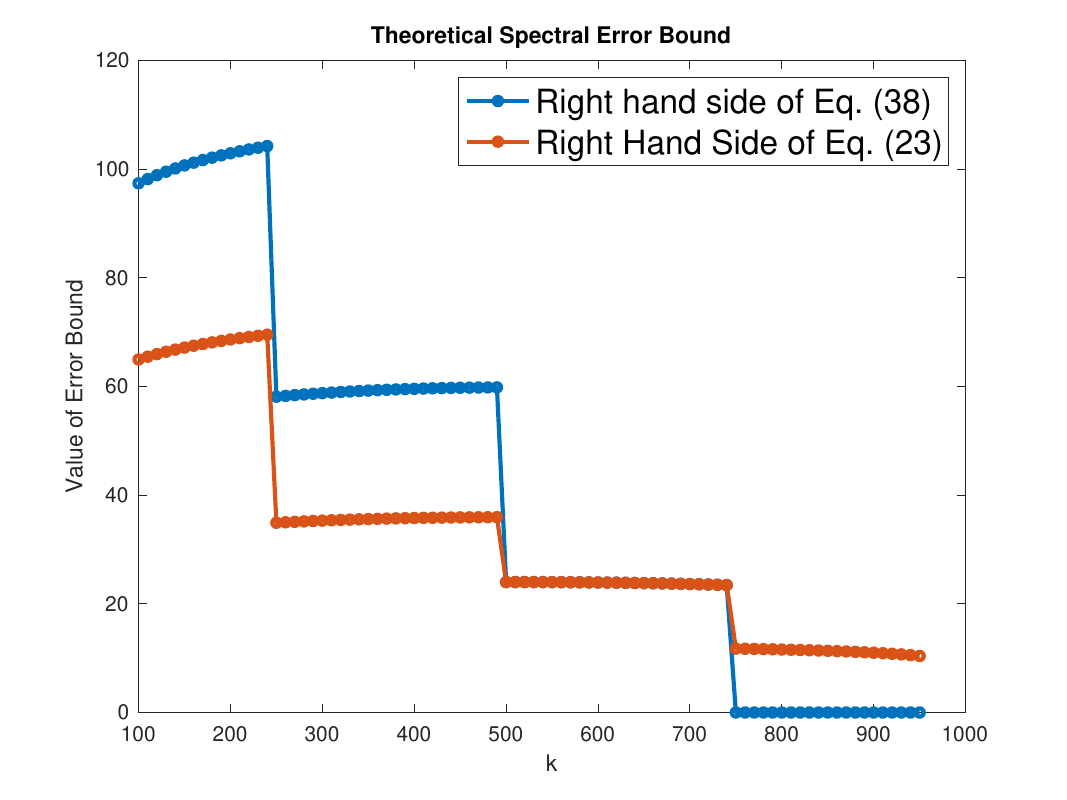}}\hfill 
		
		\vspace{-5pt}
		\caption{ \footnotesize Accuracy of projecting $X$ from Eq.~\eqref{eq: matrix where scale is better} with $n=1000$, $\beta_1=3$, $\beta_2=1$,
$\beta_3=6$, $\beta_4=2$ and $Y$ randomly generated using Algorithms~\ref{alg:approx proj} and~\ref{alg: scalled approx proj} under $l=5$, $q=2$ and varying $k \in \N$. The true PSD projections are calculated using Eq.~\eqref{eq: analytical PSD proj using eigenvalue}, where the eigen-decomposition is calculated using Matlab's \texttt{eig} function.} \label{fig: error bounds}
		\vspace{-15pt}
\end{figure}
 
\section{\textbf{PSD Projection Numerical Examples}}

 We next numerically analyze the performance of Algorithms~\ref{alg:approx proj} and~\ref{alg: scalled approx proj}, where the magnitude of the smallest eigenvalue, $\alpha$,  is computed using Alg.~\ref{alg: min eig}. We first consider the projection of $X \in \S_n$, a hand-crafted matrix defined in Eq.~\eqref{eq: matrix where scale is better}, which we showed in Corollary~\ref{cor: when sclaing bound is smaller} has a smaller projection error bound when using the scaled randomized projection (Alg.~\ref{alg: scalled approx proj}). We then benchmark our algorithms on matrices arising from practical applications, drawing from real-world sparse datasets available in the open-source University of Florida Sparse Matrix Collection~\cite{davis2011university}. The associated code and a reproducible capsule is available in~\cite{Jones_Anderson_2026}\ifsup.\else  as well as supplementary SDP numerical examples found on Arxiv in~\cite{jones2023randomizedProj}.\fi
 
 \begin{ex}
Consider the case of projecting $X \in \S_n$ from Eq.~\eqref{eq: matrix where scale is better} when $\beta_1=3, \beta_2=1$, $\beta_3=6$, $\beta_4=2$ and $n=1000$. Figure~\ref{fig: eig test} shows the eigen and singular value distribution of $X$. We note that we have chosen this test case because the spectral profile of $X$ is representative of the most ``unfriendly'' matrix class we are likely to encounter, {being full rank with no relatively small singular values. Consequently, it is expected that our proposed approximated PSD projection will exhibit reduced performance, as it involves approximating the matrix by a low-rank matrix.} Sub-figures~\ref{subfig: RHS_bound_frob_error} and~\ref{subfig: RHS_bound_spec_error} show how the error bounds in the Frobenius and spectral norms change with $k$, the number of random samples used in the Randomized Range Finder (Alg.~\ref{alg:RSNA}). These apriori error bounds are based only on the singular values and eigenvalues of $X\in \S_n$ and can be calculated without computing the full matrix projection. In both norms, we see that the red curve is initially lower than the blue curve showing that the error bound of projecting with Alg.~\ref{alg:approx proj} is smaller than that of using Alg.~\ref{alg: scalled approx proj}. This is because the first $n/4$ positive eigenvalues are equal to the leading order singular values, corresponding to $\beta_3$. However, the next $n/4$ singular values correspond to $\beta_1$ that is associated with a negative eigenvalue. For this reason it is better to use the Scaled Randomized PSD Projection (Alg.~\ref{alg: scalled approx proj}) when $k>n/4$. We see this intuition reflected in Figure~\ref{fig: error bounds} where for sufficiently large $k$ the blue curve becomes less than the red curve. This occurs in both the numerically calculated projection errors and the theoretical projection error bounds, although this is only seen for large $k$ in the theoretical projection bounds.
 \end{ex}

\begin{ex}
{We next project several matrices obtained from practical applications from the University of Florida Sparse Matrix Collection~\cite{davis2011university}. In each case $n$ exceeds 5000. As shown by the eigenvalues plotted in Fig.~\ref{fig: eigenvalues}, these matrices are low rank or approximately low rank, with many eigenvalues close to zero. The numerical performance of our projection methods is given in Table~\ref{table: ran PSD projections}. For brevity, we focus on the PSD projection error, defined as the relative Frobenius error $\|X_+ - \hat X_+\|_2 / \|X_+\|_2$, and do not consider the relative spectral error. The optimal PSD projection, $X_+$, is computed using Matlab's \texttt{eig} routine and Eq.~\eqref{eq: analytical PSD proj using eigenvalue}, whereas the RNLA projection, $\hat{X}_+$, is computed by executing either Alg.~\ref{alg:approx proj} (Vanilla) or Alg.~\ref{alg: scalled approx proj} (Scaled). Each numerical experiment represents a single reproducible run with a fixed randomization seed and parameters $l=10$ and $q=4$ (code and reproducible capsule available in~\cite{Jones_Anderson_2026}). For scaled projections, the scaling factor was set to $\alpha=|\min_i \lambda_i(X)|$ where $\alpha$ was estimated using Alg.~\ref{alg: min eig} with $N=10$. }

As expected, the projection error decreases as $k$ increases (although we maintain $k\ll n$). The scaled randomized projection, Alg.~\ref{alg: scalled approx proj}, mostly outperforms the vanilla randomized projection, Alg.~\ref{alg:approx proj}, even though the error $\|X-QQ^\top X\|_2 / \|X\|_2$ is often larger for the scaled projection;  {providing evidence our scaled PSD projection algorithm is not attempting to form a low rank approximation of the original matrix before projecting}. The \texttt{G67.mat} matrix~\cite{davis2011university} exemplifies the advantage of the scaled randomized projection, Alg.~\ref{alg: scalled approx proj}. It can be seen in Fig.~\ref{subfig: G67} that the largest singular values of \texttt{G67.mat} are equally distributed between positive and negative eigenvalues. Alg.~\ref{alg:approx proj} offers worse performance when compared to Alg.~\ref{alg: scalled approx proj} for this matrix due to the algorithm expending significant ``approximation energy" on the subspace corresponding to large negative eigenvalues. Since these components are discarded in the final PSD projection, this effectively wastes the fixed rank budget $k$ on information that does not contribute to the target solution $X_+$, whereas the scaling in Alg.~\ref{alg: scalled approx proj} redirects this ``effort" towards the relevant positive eigenspace.

 However, shifting and scaling the eigenvalue spectrum is not always advantageous. As seen with \texttt{cyl6.mat} in Fig.~\ref{subfig: cyl6} (and Table~\ref{table: ran PSD projections}), the vanilla Alg.~\ref{alg:approx proj} achieves near machine precision ($10^{-10}$) at $k=[0.25n]$, drastically outperforming the scaled version. This is  because the positive eigenvalues of \texttt{cyl6.mat} (magnitude $\approx 300$) dominate the negative ones (magnitude $\approx 100$), allowing the vanilla RNLA low-rank approximation to capture the positive eigen-subspace.

We omit empirical computation time comparisons, noting that while the theoretical $O(kn^2)$ complexity of the proposed PSD approximations suggests significant asymptotic advantages over exact methods, realizing this theoretical advantages requires optimized LAPACK implementations beyond the scope of our MATLAB-based examples. Consequently, we restrict our analysis to the reconstruction errors, which remain consistent across different code implementations.  {{For readers interested in empirical compute times, this data as well as the associated coded implementation is available in a reproducible capsule~\cite{Jones_Anderson_2026}. We note that the timings therein reflect a comparison between a sparse matrix multiplication implementation of Algorithms~\ref{alg:approx proj} and~\ref{alg: scalled approx proj} and MATLAB's dense \texttt{eig} solver.}}

\begin{table}[]   
 \scalebox{0.84}{

	{\begin{tabular}{@{} l r c 
		r r 
		r r @{}}
	\toprule
	Dataset            &    \(n\)   & \([k/n]\) 
	& \multicolumn{2}{c}{Alg.\,\ref{alg:approx proj}} 
	& \multicolumn{2}{c}{Alg.\,\ref{alg: scalled approx proj}} \\
	&            &          
	& \(\frac{\|X - QQ^\top X\|_2}{ \|X\|_2} \) 
	& \(\frac{\|X_+ - \hat X_+\|_2}{\|X_+\|_2}\) 
	& \(\frac{\|X - QQ^\top X\|_2}{\|X\|_2}\) 
	& \(\frac{\|X_+ - \hat X_+\|_2}{\|X_+\|_2}\) \\
	\midrule
	\texttt{skirt.mat}          & 12\,598     & 0.01   & 0.88   & 0.86   & 0.88    & 0.87    \\
	&             & 0.25   & 0.24    &  0.18   & 0.35    &  0.09    \\
	&             & 0.50   &  0.07    &  0.05   & 0.34    &   0.03   \\
	\addlinespace
	\texttt{shuttle\_eddy.mat}  & 10\,429     & 0.01   & 0.93   & 0.92   & 0.93   & 0.92    \\
	&             & 0.25   & 0.34   &  0.23   & 0.36   &  0.16    \\
	&             & 0.50   &  0.13   &  0.10   &  0.32   &   0.02    \\
	\addlinespace
	\texttt{cyl6.mat}           & 13\,681     & 0.01   & 0.75 & 0.72 & 0.76 & 0.72 \\
	&             & 0.25   & \(2.1\times10^{-10}\)   & \(2.3\times10^{-10}\)    & 0.33 & 0.02    \\
	&             & 0.50   & \(2.1\times10^{-15}\)   & \(1.1\times10^{-14}\)    &  0.32 & 0.01    \\
	\addlinespace
	\texttt{lowThrust\_13.mat}  & 18\,476     & 0.01   & 0.94   & 0.97   & 0.98   & 0.97   \\
	&             & 0.25   & 0.36   & 0.37   & 0.80   & 0.65   \\
	&             & 0.50   & 0.16   & 0.16   & 0.71   & 0.44   \\
	\addlinespace
	\texttt{G57.mat}           &  5\,000     & 0.01   & 0.99   & 1.00   & 0.99   &  0.97   \\
	&             & 0.25   &  0.67   &  0.71   & 0.76   &  0.41   \\
	&             & 0.50   &  0.37   &  0.40   &  0.71   &   0.04   \\
	\addlinespace
	\texttt{G67.mat}            & 10\,000     & 0.01   & 0.98   & 0.99   & 0.99   & 0.97   \\
	&             & 0.25   & 0.67   & 0.71   & 0.76   &  0.41   \\
	&             & 0.50   &  0.37   &  0.39   & 0.71   &   0.04   \\
	\bottomrule
\end{tabular}}} \vspace{-4pt}
		\caption{\footnotesize {Errors} associated with finding RNLA and deterministic PSD projections of benchmark matrices from~\cite{davis2011university} {(code and reproducible capsule available in~\cite{Jones_Anderson_2026})}.} 	\vspace{-10pt}
		\label{table: ran PSD projections}
\end{table} 
\end{ex}

\ifthenelse{\boolean{longver}}{	\begin{figure*}
		%
		\subfloat[ \texttt{\scriptsize  shuttle$\_$eddy.mat} \label{subfig: shuttle eddy}]{\includegraphics[width=0.49 \linewidth, trim = {1cm 0cm 1cm 0cm}, clip]{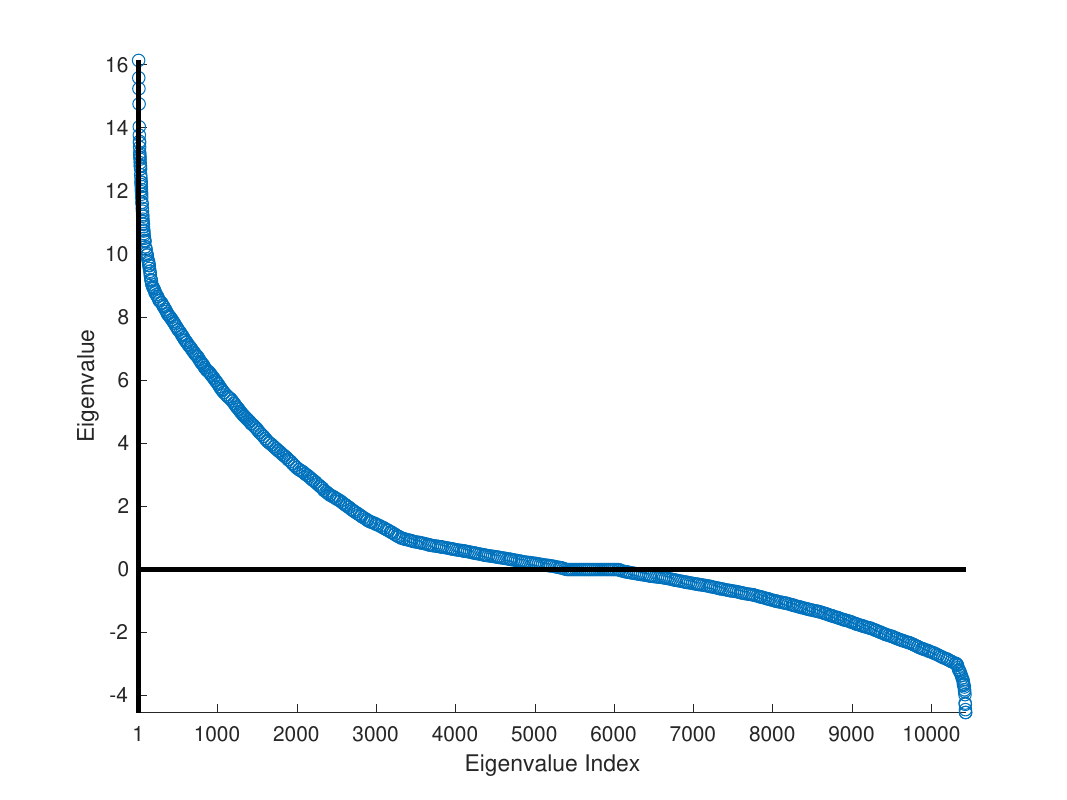}}\hfill
		\subfloat[ \texttt{\scriptsize cyl6.mat} \label{subfig: cyl6}]{\includegraphics[width=0.49 \linewidth, trim ={1cm 0cm 1cm 0cm}, clip]{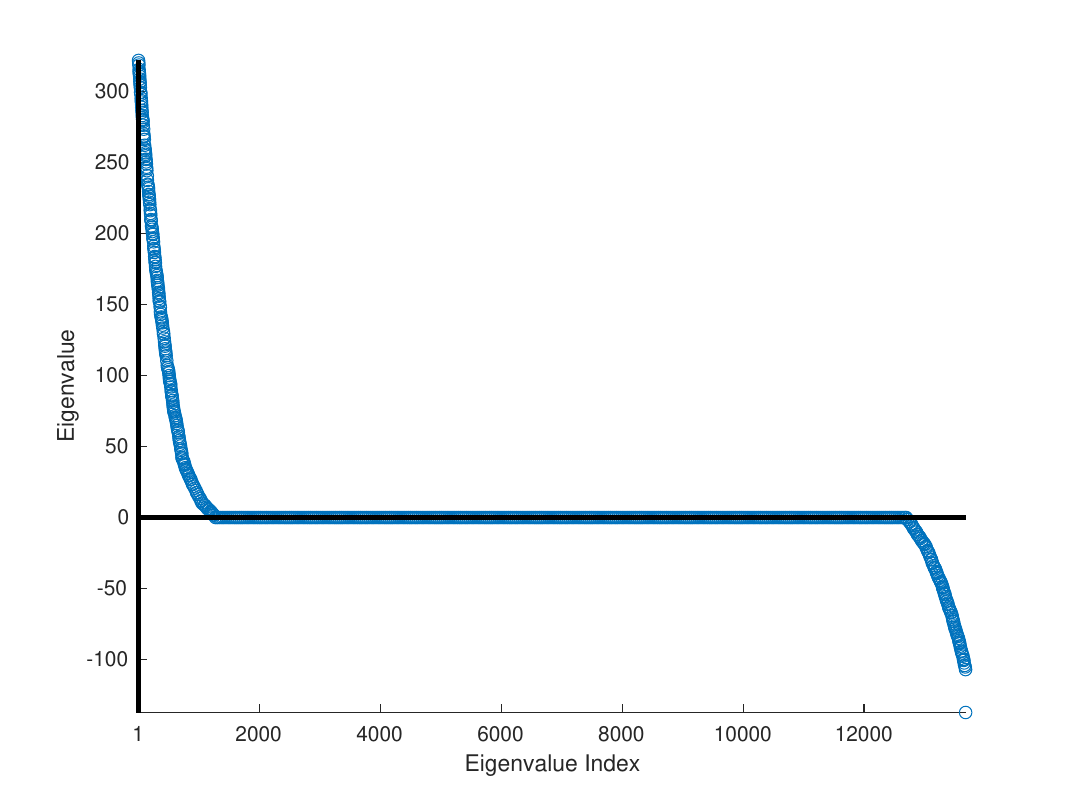}}\hfill 
\subfloat[ \texttt{\scriptsize lowThrust\_13.mat} \label{subfig: low thrust plot} ]{\includegraphics[width=0.49 \linewidth, trim = {1.5cm 0cm 1cm 0cm}, clip]{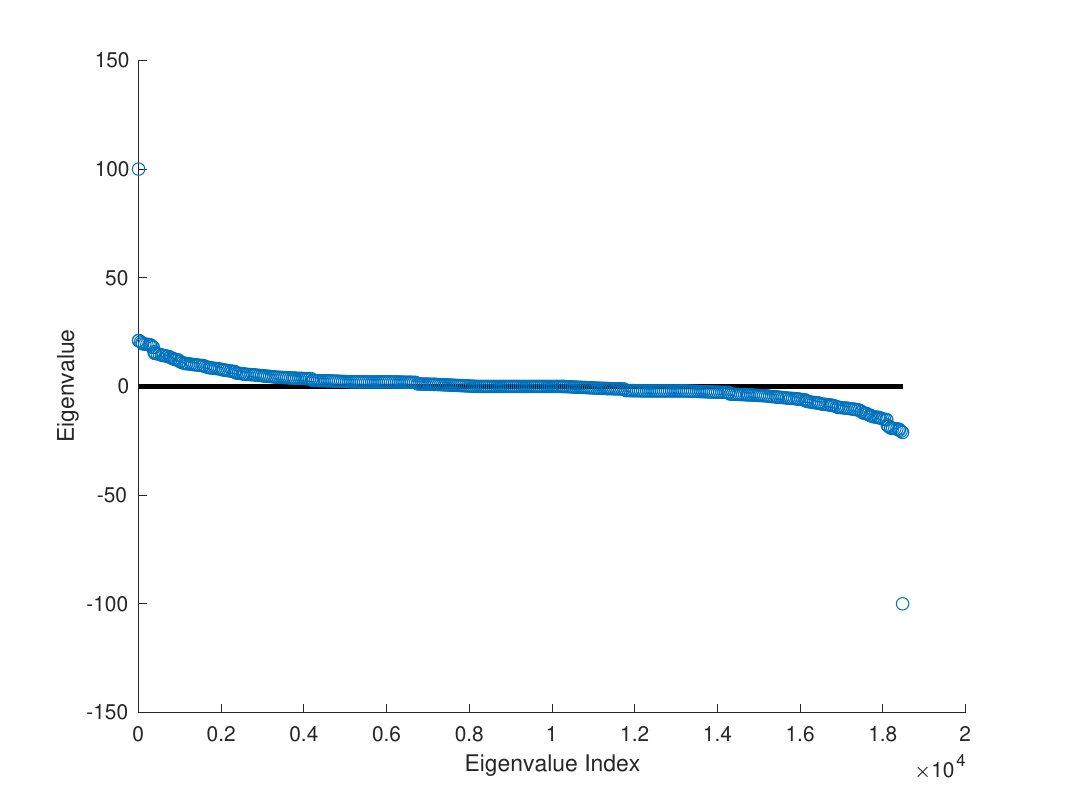}} \hfill 
		\subfloat[ \texttt{\scriptsize G67.mat} \label{subfig: G67}]{\includegraphics[width=0.49 \linewidth, trim = {1cm 0cm 1cm 0cm}, clip]{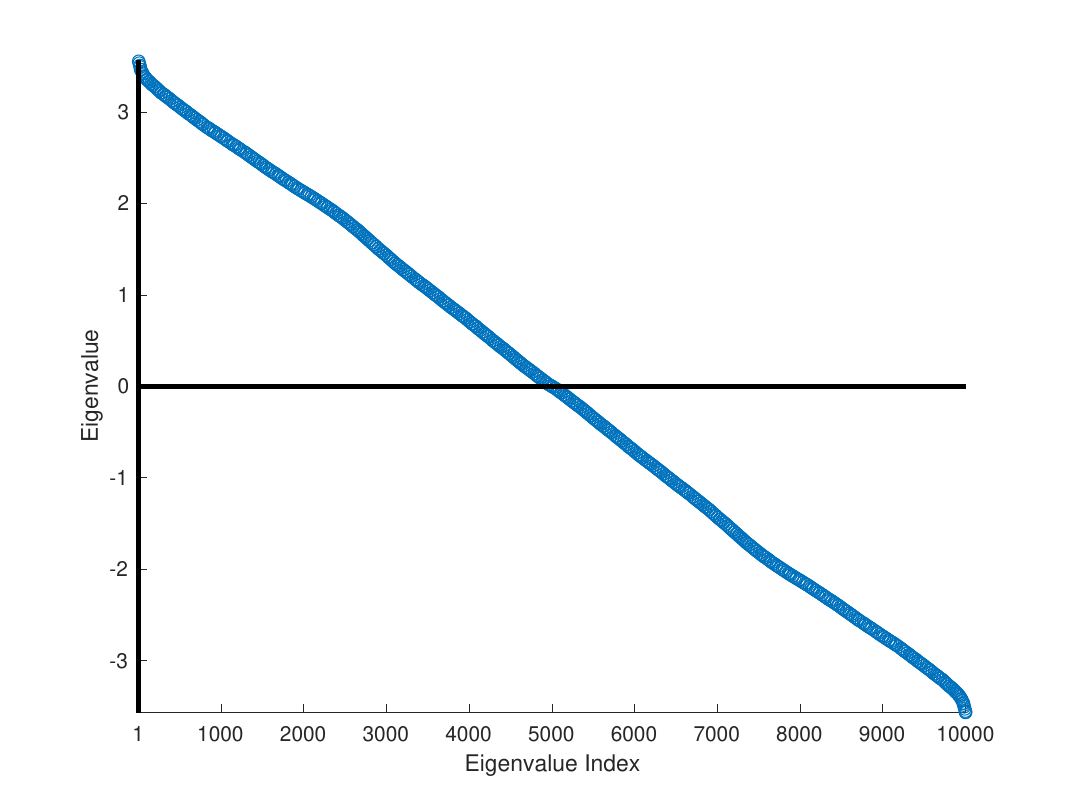}}\hfill
		\caption{\footnotesize Scatter plot showing the distribution of the eigenvalues $\lambda_1 \ge \dots \ge \lambda_n$ for various benchmark matrices from~\cite{davis2011university}, whose projections are analysed in Table~\ref{table: ran PSD projections}. {The error bounds in Corollary~\ref{cor: error bound on rand projection} and Proposition~\ref{prop: proj error 2} depend on the eigenvalues and singular values of the matrices. These plots visualize the spectral properties of a sample of our benchmark matrices, allowing us to anticipate how well they can be approximated by a low-rank structure. By examining the distribution of eigenvalues, we can use the derived error bounds to predict which approximate projection method (scaled or non-scaled) is likely to perform better. The true numerical errors given in Table~\ref{table: ran PSD projections}.}
			} \label{fig: eigenvalues}
		\vspace{-5pt}
\end{figure*} }

		\vspace{-0.3cm}
		\section{Conclusion}
		\vspace{-0.1cm}
		We have proposed a highly scalable method of harnessing RNLA to compute approximate PSD projections. We {have} shown that in some cases it is advantageous to scale the matrix to align the largest singular values with the largest positive eigenvalues before approximating by a low rank matrix and projecting. {We have provided several probabilistic error bounds for our proposed approximate PSD projection methods, which demonstrate that increasing the approximation rank $k$ leads to more accurate projections. The choice of $k$ is therefore critical in balancing computational efficiency with projection quality. Future works involves exploring replacing Alg.~\ref{alg:RSNA} with the adaptive rank revealing algorithms of~\cite{martinsson2016randomized, yu2018efficient} and investigating integration into SDP first order solvers.}
		


  \section*{Acknowledgement}
The first author would like to thank Dr Jeff Wadsworth – Battelle Knowledge Exchange Scheme Award, the University of Sheffield and Columbia University for supporting a summer visit that resulted in a collaboration and the production of this work. The second author would like to acknowledge funding from the National Science Foundation under grant ECCS 2144634.
		\ifthenelse{\boolean{longver}}{\vspace{-0.2cm}}{
			\vspace{-0.1cm}}
		\bibliographystyle{ieeetr}
		\bibliography{bib_Randomized_PSD_projections}
		\ifthenelse{\boolean{longver}}{ 
			\section{Appendix}

			\begin{lem}[{Frobenius norm matrix absolute value inequality}] \label{lem: An Inequality for the Frobenius norm}
				Suppose $A,B \in \S_n$ then
				\vspace{-0.2cm}\begin{align} \label{ineq: An Inequality for the Frobenius norm}
					\left|\left| \left(A^\top A \right)^{\frac{1}{2}}   - \left(B^\top B \right)^{\frac{1}{2}} \right| \right|_2 \le \| A -B\|_2.
				\end{align}
			\end{lem}\begin{proof}
				In finite dimensions the Hilbert-Schmidt operator is identical to the Frobenius norm and hence Ineq.~\eqref{ineq: An Inequality for the Frobenius norm} follows by Lemma 5.2 in~\cite{araki1971quasifree}. {Alternatively, this inequality appears in Corollary VII.5.6 in~\cite{bhatia2013matrix}.}
			\end{proof}
			
			\begin{lem}[{Spectral norm matrix absolute value inequality}] \label{lem: An Inequality for the spectral norm}
				Suppose $A,B \in \S_n$ then
				\vspace{-0.2cm}\begin{align} \label{ineq: An Inequality for the spectral norm}
					& \left|\left| \left(A^\top A \right)^{\frac{1}{2}}   - \left(B^\top B \right)^{\frac{1}{2}} \right| \right|_\infty\le \frac{2}{\pi} \left( 2 + \log \left( \frac{\|A\|_\infty +\|B\|_\infty}{\|A-B\|_\infty} \right) \right) \| A -B\|_\infty.
				\end{align}
			\end{lem}
		\begin{proof}
		{This inequality first appears in \cite{kato1973continuity} but is also present in Theorem X.2.5~\cite{bhatia2013matrix}.}
		\end{proof}
			
			


			
			
			\begin{thm}[\cite{halko2011finding}] \label{thm: expected error bounds of proto alg}
				Consider $A\in \R^{m \times n}$. Choose a target rank $k  \ge  2$ and an oversampling parameter $l \ge  2$, where $k + l \le  \min\{m, n\}$. Execute the {Randomized Range Finder} (Alg.~\ref{alg:RSNA}) for $q \in \N$ to output an orthonormal approximate basis $Q$ which satisfies,
				\vspace{-0.2cm}\begin{align} \label{ineq: proto error bound in spectral radius}
					&\mathbb{E}[\|A-Q Q^\top A\|_\infty ]\\ \nonumber 
					&\le \hspace{-0.1cm}  \left(1+\sqrt{\frac{k}{l-1}} + \frac{e \sqrt{k+l}}{l}\sqrt{\min\{m,n\}-k} \right)^{\frac{1}{2q+1}} \hspace{-0.4cm} \sigma_{k+1}(A),
				\end{align}
				where the expectation is taken with respect to the randomly generated matrix $\Omega$.
				
				Moreover when the {Randomized Range Finder} (Alg.~\ref{alg:RSNA}) is executed with $q=0$ then the following bound in the Frobenius norm is satisfied,
				\begin{align} \label{ineq: randomized frobenious error bound}
					&\mathbb{E}[\|A-Q Q^\top A\|_2 ]\le \sqrt{\left(1+\frac{k}{l-1} \right) \left(\sum_{j>k} \sigma_j(A)^2  \right)}.
				\end{align}
		\end{thm} }

\begin{lem}
For any $x\ge 0$ we have
    \begin{align}   \label{ineq: logx bound}
    \log(x) & \le \frac{1}{2} \sqrt{x} +1 \quad \text{and}\\
    \label{ineq: xlogx bound}
        -x\log(x) &\le \min\{ \sqrt{x},e^{-1}\},
    \end{align}
    where $\log$ is the natural logarithm. 
\end{lem}
\begin{proof}
Let $g(x)=\log(x)-0.5 \sqrt{x}-1$. To show Ineq.~\eqref{ineq: logx bound} we show $\sup_{x \ge 0} g(x) \le 0$. To find the max point of $g$ we find it's stationary points, $g'(x)=\frac{1}{x}-\frac{1}{4\sqrt{x}}=\frac{4-\sqrt{x}}{4 x}$. Clearly $g'(16)=0$. Moreover, $g''(x)=-\frac{1}{x^2}+\frac{1}{8x^{3/2}}=\frac{\sqrt{x}-8}{8x^2}$ and hence $g''(16)=\frac{-4}{8(16)^2}<0$ implying this stationary point is a max point. Therefore, $\sup_{x \ge 0} g(x)=g(16)=\log(16)-2-1=-0.227<0$.

    Let $f(x):=-x\log(x) - \sqrt{x}$. To show Ineq.~\eqref{ineq: xlogx bound} we show $f(0)=0$ and $f'(x)<0$, that is $f$ is monotonically decreasing and hence $0=f(0) \ge f(x)$. To show $f(0)=0$ note that $\sqrt{0}=0$ and $x \log(x)$ is continuous with $\lim_{x \to 0} \{x \log(x)\} =\lim_{x \to 0} \left\{ \frac{\log(x)}{1/x} \right\}=0$ can be shown by L'Hôpital's rule. Now, $f'(x)= -\log(x) -1 - \frac{1}{2 \sqrt{x}}$. Clearly, when $x>1$ we have $\log(x)>0$ and hence each of the terms in $f'(x)$ are negative implying $f'(x)$ is negative when $x>1$. We next show $f'(x)<0$ when $x \in [0,1]$ which is equivalent to showing $f'(1/y)=\log(y) -1 -(1/2)\sqrt{y}$ is negative for $y \ge 1$. By Ineq.~\eqref{ineq: logx bound} we have that $f'(1/y)=\log(y) - (1/2)\sqrt{y} -1 \le -1<0$ and hence $-x\log(x)<\sqrt{x}$. On the other hand the maximum of $h(x):=-x\log(x)$ occurs when $h'(x)=-\log(x) -1 =0$ which is solved at $x=1/e$. Moreover, $h''(1/e)=-\frac{1}{1/e}=-e<0$ implying this is a max point. Hence $-x\log(x)<h(1/e)=e^{-1}$.  
\end{proof}

\begin{lem} \label{lem: interchange PSD proj}
    Consider any $X \in \S_k$. Suppose $Q \in \R^{n \times k}$ is such that $Q^\top Q=I$. Then
    \begin{align} \label{eq: interchange PSD projection}
        QX_+Q^\top =(QXQ^\top)_+
    \end{align}
\end{lem}
\begin{proof}
    \begin{align} \label{pfeq: 11}
        QX_+Q^\top \marksymb{=}{\ref{eq: projection in polar form}} \frac{1}{2}Q(X + \sqrt{X^\top X} ) Q^\top=\frac{1}{2}(QXQ^\top + Q\sqrt{X^\top X}Q^\top).
    \end{align}
    We next show $Q\sqrt{X^\top X}Q^\top= \sqrt{QX^\top XQ^\top}$. First note that $Q\sqrt{X^\top X}Q^\top \in S_n^+$ by definition that $\sqrt{X^\top X} \in \S_n^+$ is the uniquely defined square root of $X^\top X$ (see notation section). Next, note that it follows from $Q^\top Q=I$ that \[Q\sqrt{X^\top X}Q^\top Q\sqrt{X^\top X}Q^\top=Q\sqrt{X^\top X}\sqrt{X^\top X}Q^\top=Q X^\top X Q .\] 
    Since the square root of $Q X^\top X Q$ is uniquely defined it follows that $\sqrt{QX^\top XQ^\top}=Q\sqrt{X^\top X}Q^\top$. Upon substituting into Eq.~\eqref{pfeq: 11} and using $Q^\top Q=I$ again, we deduce that
        \begin{align*}
        QX_+Q^\top & \marksymb{=}{\ref{pfeq: 11}}\frac{1}{2}(QXQ^\top + Q\sqrt{X^\top X}Q^\top)=\frac{1}{2}(QXQ^\top + \sqrt{QX^\top XQ^\top})\\
        &=  \frac{1}{2}(QXQ^\top + \sqrt{QX^\top Q^\top Q XQ^\top}) \marksymb{=}{\ref{eq: projection in polar form}} (QXQ^\top)_+.
    \end{align*}
\end{proof}

\begin{lem} \label{lem: spectal projection bound}
 {	Consider $A,B \in \S_n$. Suppose $Q \in \R^{n \times k}$ is such that $Q^\top Q=I$. Then
	\begin{align} \label{eq: inf norm bound}
		\|QQ^\top A QQ^\top  + (I-QQ^\top )B (I-QQ^\top ) \|_{\infty} \le \max\{\|A\|_{\infty},\|B\|_\infty  \}.
	\end{align} }
\end{lem}
\begin{proof}
 {	Let $X := QQ^\top A QQ^\top + (I-QQ^\top )B (I-QQ^\top )$. Since $A, B \in \S_n$, $X$ is also symmetric. Thus, the spectral norm is equivalent to the induced $l_2$ operator norm $\|X\|_{\infty} = \sup_{ v \ne 0} \left\{ \frac{|v^\top X v|}{\|v\|_2^2} \right \}$ (see Chapter 5.6 in~\cite{horn2012matrix}).
	
	For any $v \in \R^n$, let $v_1 := QQ^\top v$ and $v_2 := (I-QQ^\top) v$. It is clear that $v = v_1 + v_2$, and since $QQ^\top (I-QQ^\top) = 0$ we have $v_1 ^\top v_2=0$. It then follows that 
	\begin{align} \label{pfeq:pythag}
\|v\|_2^2 = \|v_1\|_2^2 +2v_1^\top v_2+ \|v_2\|_2^2=\|v_1\|_2^2 + \|v_2\|_2^2.
	\end{align}
	
	Substituting $v = v_1 + v_2$ into the quadratic form, the cross terms vanish because $QQ^\top (I-QQ^\top) = 0$ and using the triangle inequality we get that,
	\begin{align*}
		|v^\top X v| &= |(v_1 + v_2)^\top X (v_1 + v_2)| \\
		&= |v_1^\top QQ^\top A QQ^\top  v_1 + v_2^\top (I-QQ^\top )B (I-QQ^\top ) v_2| \\
		&\le |v_1^\top QQ^\top A QQ^\top  v_1| + |v_2^\top (I-QQ^\top )B(I-QQ^\top ) v_2|.
	\end{align*}
	
	By the definition of the spectral norm for symmetric matrices, we have
	\begin{align*}
	|v_1^\top QQ^\top A QQ^\top v_1| & \le \|QQ^\top A QQ^\top \|_\infty \|v_1\|_2^2 \le \|QQ^\top\|_{\infty}^2 \|A\|_\infty \|v_1 \|_2^2\\
	& \marksymb{=}{\ref{eq: QQ 2 norm} }\|A\|_\infty \|v_1 \|_2^2.
	\end{align*}  and  similarly,
	\begin{align*}
|v_2^\top (I-QQ^\top )B(I-QQ^\top ) v_2| \le \|B\|_\infty \|v_2\|_2^2
	\end{align*}

Substituting these in:
	\begin{align*}
		|v^\top X v| &\le \|A\|_\infty \|v_1\|_2^2 + \|B\|_\infty \|v_2\|_2^2 \\
		&\le \max\{\|A\|_\infty, \|B\|_\infty\} \left( \|v_1\|_2^2 + \|v_2\|_2^2 \right) \\
		&\marksymb{=}{\ref{pfeq:pythag}} \max\{\|A\|_\infty, \|B\|_\infty\} \|v\|_2^2.
	\end{align*}
	
	Dividing by $\|v\|_2^2$ (for $v \ne 0$) and taking the supremum over all $v$, we get:
	\begin{align*}
		\|X\|_{\infty} = \sup_{v \ne 0} \frac{|v^\top X v|}{\|v\|_2^2} \le \max\{\|A\|_{\infty}, \|B\|_\infty \}.
	\end{align*} }
\end{proof}

\begin{lem} \label{lem: caluclate min eigenvalue}
    Consider $X \in \S_n$. It follows that $|\lambda_n(X)|:=|\min_i \lambda_i(X)|=|\sigma_{1}(Y)-\sigma_{1}(X)|$, where $Y:=X-\sigma_{1}(X)I$.
\end{lem}
\begin{proof}
Since $X \in \S_n$ is a symmetric matrix it has the property that the absolute values of their eigenvalues are exactly equal to their singular values, i.e., $\{|\lambda_1(X)|,\dots,|\lambda_n(X)| \}=\{\sigma_1(X),\dots,\sigma_n(X)\}$. It follows that $Y \in \S_n$ is simultaneously diagonalizable with $X$ with eigenvalues shifted by negative $\sigma_1(X)$. That is, $\lambda_i(Y)= \lambda_i(X)-\sigma_1(X)=\lambda_i(X)-\max_i |\lambda_i(X)| \le 0$ and hence $\sigma_i(Y)=|\lambda_{n-i+1}(X)-\sigma_1(X)|$. There are two possible cases depending if the largest singular value of $X$ corresponds to the absolute value of the maximal or minimal eigenvalue of $X$:
    \begin{itemize}
		\item \textbf{\underline{Case 1:}}  $\sigma_{1}(X)=| \lambda_n(X)|$
  
  This implies that the absolute value of the minimal eigenvalue is greater than the absolute value of the maximum eigenvalue of $X$, $|\lambda_n(X)| \ge |\lambda_{1}(X)|$. Then it follows that the minimal eigenvalue of $X$ must be negative, $\lambda_n(X)\le0$, otherwise if $0\le\lambda_n(X) \le \lambda_{1}(X)$ we would contradict $|\lambda_n(X)| \ge |\lambda_{1}(X)|$. Now $\sigma_{1}(Y)=| \lambda_n(X) -\sigma_1(X)|=| \lambda_n(X) -|\lambda_n(X)\|=2|\lambda_n(X)|$. Hence $|\sigma_{1}(Y)-\sigma_{1}(X)|=| (2|\lambda_n(X)|- | \lambda_n(X)|)|=| \lambda_n(X)|$.
		\item \textbf{\underline{Case 2:}}  $\sigma_{1}(X)=|\lambda_{1}(X)|$ 
  
  By a similar argument to Case 1 we get $\lambda_{1}(X)\ge0$ and hence $\sigma_1(X)=\lambda_{1}(X)$. Moreover since $\lambda_{1}(X) \ge  \lambda_n(X)$ we get that $\sigma_{1}(Y)=\lambda_{1}(X)- \lambda_n(X)$. Hence, $|\sigma_{1}(Y)-\sigma_{1}(X)|=|\lambda_{1}(X)- \lambda_n(X) - \sigma_{1}(X)|=|\lambda_n(X)|$.
	\end{itemize}
\end{proof}

				\ifsup
				\section{Supplementary Material: SDP Numerical Examples}
				In Section~\ref{sec: Projecting onto the Cone of Positive Semidefinite Matrices} we introduced two randomized algorithms for approximate projecting matrices onto the PSD cone (Algorithms~\ref{alg:approx proj} and~\ref{alg: scalled approx proj}). PSD projections are widely used in many first-order SDP solvers. {To demonstrate the ease of integrating our approximate PSD projection methods into existing first-order SDP algorithms, we incorporate them into a specific solver, which we describe in detail in the following section, and conduct numerical experiments to assess their impact on solution quality.  }

				 We now recap the first-order solver for solving Semidefinite Programming Least Square{s} (SDLS) problems from~\cite{malick2004dual}, showing how to modify this algorithm to solve {regularized} SDPs and how we can employ our randomized PSD projection algorithms. Consider the following SDLS problem,
				\ifthenelse{\boolean{longver}}{\begin{align} \label{opt: SDLS}
								&\min_{X \in \S_n^+} \frac{1}{2} \left| \left|X - \frac{1}{\rho} C \right| \right|_2^2\\ \nonumber
								\text{subject to: } & \Tr(A_i^\top X)=b_i \text{ for } i \in \{1,\dots,m\},
						\end{align} }{\vspace{-0.15cm}  \begin{align} \label{opt: SDLS}
								&\min_{X \in \S_n^+} \frac{1}{2} \left| \left|X - \frac{1}{\rho}  C \right| \right|_2^2\\ \nonumber
								\text{subject to: } & \Tr(A_i^\top X)=b_i \text{ for } i \in \{1,\dots,m\},
						\end{align} }
						where $\rho>0$, $b \in \R^m$, $C \in \S_n$ and $\{A_i\}_{i \in \{1,\dots,m\}} \subset \S_n$ {is the problem data, and $X$ the decision variable}. 
						
						
						{Following  the partial dualization approach of \cite{malick2004dual,malick2009regularization} and~\cite{henrion2011projection} for solving~\eqref{opt: SDLS}}. Consider the {partial} Lagrangian of Opt.~\eqref{opt: SDLS}:
						\ifthenelse{\boolean{longver}}{\begin{align*} 
										&	L(X;y) :=\frac{1}{2}  \left| \left|X - \frac{1}{\rho}  C \right| \right|_2^2 - \sum_{i=1}^m y_i \left(\Tr(A_i^\top X) - b_i \right),
								\end{align*}}{
									\vspace{-0.35cm} \begin{align*} 
											&	L(X;y) :=\frac{1}{2}  \left| \left|X - \frac{1}{\rho}   C \right| \right|_F^2 - \sum_{i=1}^m y_i \left(\Tr(A_i^\top X) - b_i \right),
									\end{align*}}
									where only the affine constraint{s} have been pulled into the Lagrangian. 
											The dual problem is thus given by,
											\ifthenelse{\boolean{longver}}{}{
													\vspace{-0.1cm}}\begin{align} \label{opt: dual problem} 
													\sup_{y \in \R^m}  \inf_{X \in \S_n^+} L(X;y).
												\end{align}

											Define $X^*(y):=\inf_{X \in \S_n^+} L(X;y)$. The following lemma gives an analytical expression for $X^*(y)$.
											\begin{lem} \label{lem: analytical sol to dual}
													Consider $X^*(y):=\inf_{X \in \S_n^+} L(X;y)$. It follows that
													\ifthenelse{\boolean{longver}}{\begin{align} \label{eq: analytical sol. of dual}
																	X^*(y)= \left(\frac{1}{\rho} C + \sum_{i=1}^my_i A_i \right)_+.
															\end{align} }{
															\vspace{-0.4cm} \begin{align} \label{eq: analytical sol. of dual}
																	X^*(y)= \left(\frac{1}{\rho}  C + \sum_{i=1}^my_i A_i \right)_+.
															\end{align} }
												\end{lem}
											\begin{proof} 
													\ifthenelse{\boolean{longver}}{%
															By completing the square in $L(x;y)$ it follows that,
															\begin{align*}
																	L(X;y) \hspace{-0.05cm}= \hspace{-0.05cm} \frac{1}{2} & \left| \left|X - \frac{1}{\rho} C - \sum_{i=1}^my_i A_i \right| \right|_2^2 \hspace{-0.2cm} - y^\top b  - \frac{1}{2}\left| \left|\sum_{i=1}^my_i A_i \right| \right|_2^2 \hspace{-0.2cm}- \sum_{i=1}^m \frac{y_i}{\rho} \Tr(A_i^\top C).
																\end{align*}
															Only the first term in $L$ involves $X$. Hence, solving $\inf_{X \in \S_n^+} L(X,y)$ is equivalent to solving $\inf_{X \in \S_n^+} \left| \left|X - \frac{1}{\rho} C - \sum_{i=1}^my_i A_i \right| \right|_2^2$. Therefore by Theorem~\ref{thm: Analytical PSD Projection} Eq.~\eqref{eq: analytical sol. of dual} holds. {Furthermore, note that as $A_i$ and $C$ are symmetric, so is their weighted sum, and hence the projection in solution~\eqref{eq: analytical sol. of dual} exists. }
														}{%
															{   Complete the square of $L$ and apply Thm.~\ref{thm: Analytical PSD Projection}. }
														}
												\end{proof}

											Now by Lemma~\ref{lem: analytical sol to dual} the dual problem, given in Eq.~\eqref{opt: dual problem}, is reduced to the following unconstrained optimization problem,
											\ifthenelse{\boolean{longver}}{}{
													\vspace{-0.45cm}} \begin{align} \label{opt: unconstrained opt for sdls}
													\sup_{y \in \R^m} \theta(y):= L(X^*(y),y).
												\end{align}
											It has been previously shown in~\cite{malick2004dual} that Opt.~\eqref{opt: unconstrained opt for sdls} is concave, coercive, and has Lipschitz continuous gradient that is given by 
											\ifthenelse{\boolean{longver}}{}{
													\vspace{-0.1cm}}
										    \begin{align*}
											    \nabla \theta(y)=-(\Tr(A_1^\top X^*(y)) - b_1,\dots, \Tr(A_m^\top X^*(y)) - b_m )^\top. 
											    \end{align*}
													Using a first order method, like Gradient Descent (GD), to solve Opt.~\eqref{opt: unconstrained opt for sdls} necessitates the computation of PSD projections through $X^*$ in Eq.~\eqref{eq: analytical sol. of dual}. For large scale problems this may be numerically prohibitive. Therefore we next propose computing $\nabla \theta$ by using our proposed randomized projection from Sec.~\ref{sec: Projecting onto the Cone of Positive Semidefinite Matrices}. This approach is illustrated in Alg.~\ref{alg:Dual Gradient Descent}. The expected error of the computation of $\nabla \theta$ through our RNLA PSD projection is given in the following corollary.
													\ifthenelse{\boolean{longver}}{}{
															\vspace{-0.1cm}}
													\begin{cor}[Randomized gradient error bound] \label{cor: approx dual grad}
															Let $y \in \R^m$ and consider $X:=\left(\frac{1}{\rho} C + \sum_{i=1}^my_i A_i \right) \in \S_n$, where $C, A_i \in \S_n$, and $\rho>0$. Choose a target rank $k  \ge  2$ and an oversampling parameter $l \ge  2$, where $k + l \le n$. Execute Alg.~\ref{alg:approx proj} with $q=0$ to output $\hat{X}_+$. Construct an approximation of  $\nabla \theta$; 
															\ifthenelse{\boolean{longver}}{}{
																	\vspace{-0.1cm}}\begin{align*}
																	{\nabla \hat\theta(y)}:=-(\Tr(A_1^\top \hat{X}_+) - b_1,\dots, \Tr(A_m^\top  \hat{X}_+) - b_m )^\top.
																\end{align*}
															%
															\vspace{-0.6cm}
												   \text{Then} \begin{align} \label{eq: grad error bound expectation}
																		&	\mathbb{E}[\|\nabla \theta(y)  - {\nabla \hat\theta(y)}\|_2] \hspace{0cm}  \le  \sqrt{2}\left({\sum_{i=1}^m \|A_i\|_2} \right) \eps_1(\{\sigma_i(X)\}_{i=1}^n,k,l) 
																\end{align} 
															where the expectation is taken with respect to the randomly generated matrix $\Omega$ in Alg.~\ref{alg:RSNA} and the $\eps_1$ term is defined in Eqs~\eqref{eps1}. 
														\end{cor}
													\begin{proof}
												For convenience, let us denote the Frobenius inner product by $<X,Y>_F=\Tr(X^\top Y)$. Now, applying the Cauchy Schwartz inequality and the sub-addativity property of the square root operator we get that 
												\begin{align*}
												   & \|\nabla \theta(y) -\nabla \hat{\theta}(y) \|_2=\|[\Tr(A_1^\top (\hat{X}_+-X_+)),\dots,\Tr(A_1^\top (\hat{X}_+-X_+)) ]^\top \|_2\\
												    &=\sqrt{\sum_{i=1}^m \langle A_i,\hat{X}_+-X_+\rangle_F^2} \quad \le \quad \sqrt{\sum_{i=1}^m \|A_i\|_2^2 \|\hat{X}_+-X_+\|_2^2 } \\
												    &\le \sum_{i=1}^m \|A_i\|_2 \|\hat{X}_+-X_+\|_2.
												\end{align*}
												  Finally, take expectations and apply Ineq.~\eqref{ineq: F proj error}.
														\end{proof}
													
													By a similar argument to Corollary~\ref{cor: approx dual grad} we can also establish an error bound for using the randomized scaled projection (Alg.~\ref{alg: scalled approx proj}) using Proposition~\ref{prop: proj error 2}. We can also derive bounds in terms of the spectral norms for both algorithms using the equivalence between the two norms, $\|\hat{X}_+-X\|_2 \le \sqrt{n} \|\hat{X}_+-X\|_\infty$, and the spectral error bounds in Ineqs~\eqref{ineq: spectral proj error} and~\eqref{ineq: proj error 3}. These gradient error bounds are required to use the many convergence results for problems with stochastic biased gradient error~\cite{chen2018stochastic,hu2020biased,scaman2020robustness}.
													\begin{algorithm} 
															\caption{Dual gradient descent to solve Opt.~\eqref{opt: SDLS}.}\label{alg:Dual Gradient Descent}
															\hspace*{\algorithmicindent} \textbf{Input:} \text{SDP parameters:} $C \hspace{-0.05cm} \in \hspace{-0.05cm} \S_n$, $\rho \hspace{-0.05cm} > \hspace{-0.05cm}0$, $A_i \hspace{-0.05cm}\in \hspace{-0.05cm} \S_n$, $b \hspace{-0.05cm} \in \hspace{-0.05cm} \R^m$. \Comment{{Opt.~\eqref{opt: SDLS}}} \\
															\hspace*{\algorithmicindent} \text{RNLA proj parameters:} {Approx rank} $k\in \N$, {oversampling parameter} $l\in \N$,\\
															\text{ } \hspace{4.5cm} {power iteration parameter} $q \in \N$,\\ \text{ } \hspace{4.5cm} { number of power iterations} $N \in \N$, \\
																\text{ } \hspace{4.5cm} {and choice of projection algorithm} scal$ \in \hspace{-0.05cm} \{{-1},0, \hspace{-0.05cm} 1\}$.\\
															\hspace*{\algorithmicindent} \text{GD parameters:} $\eps>0$ {(stopping condition)}, $\beta>0$ {(learning rate)}\\
															\text{ } \hspace{3.1cm} {and } $M \in \N$ {(max number of iterations)}.\\
															\hspace*{\algorithmicindent} \textbf{Output:} ${\hat{X}_+}\in \S_n^+$ approximate solution to Opt.~\eqref{opt: SDLS}. 
															\begin{algorithmic}[1]
																	\State $y=\texttt{randn}(m,1)$  \Comment{Randomized GD initialization}
																	\State $\nabla \theta(y)=\eps {\times} \texttt{ones}(m,1);\text{ }  i=0$  \Comment{Gradient initialization}
																	\While{$\|\nabla\theta(y)\|_2> \eps$ and $i \le M$}
																	\State$ i=i+1$
																	\State $X=(1/\rho) C + \sum_{i=1}^my_i A_i $
																				\If{\text{scal}=-1} 
																				\State {$X=UDU^\top$} \Comment{Full eigen-decomposition}
																	\State {$\hat{X}_+=UD_+U^\top$}
																	\Comment{Eq.~\eqref{eq: analytical PSD proj using eigenvalue}: Exact projection}
																		\ElsIf{\text{scal}=0} 
																	\State $\hat{X}_+=\texttt{ran$\_$proj}(X,k,l,q)$ \Comment{Alg.~\ref{alg:approx proj}}
																	\ElsIf{\text{scal}=1}
																	\State $\alpha=\texttt{min$\_$eig}(X,N)$ \Comment{Alg.~\ref{alg: min eig}}
																	\State $\hat{X}_+=\texttt{ran$\_$proj$\_$scal}(X,k,l,q,\alpha)$ \Comment{Alg.\ref{alg: scalled approx proj}}
																	\EndIf
																	\State  $\nabla \hspace{-0.05cm} \hat \theta(y) \hspace{-0.05cm}= \hspace{-0.05cm}(\Tr(A_1^\top \hat{X}_+) , \hspace{-0.05cm} \dots \hspace{-0.05cm}, \hspace{-0.1cm} \Tr(A_m^\top \hat{X}_+) )^\top \hspace{-0.15cm}-\hspace{-0.05cm} b$  \Comment{Cor.~\ref{cor: approx dual grad} }
																	\State $y = y - \beta \nabla \theta(y)$  \Comment{Update GD iteration}
																	\EndWhile
																\end{algorithmic}
														\end{algorithm}
													{We now reformulate Opt.~\eqref{opt: SDLS} into a more natural form for specifying SDPs.}
													By expanding the squared objective, multiplying through by $\rho>0$ and removing constant terms in the SDLS problem given in Opt.~\eqref{opt: SDLS} we arrive at the following equivalent problem,
													\ifthenelse{\boolean{longver}}{}{
															\vspace{-0.1cm}} \begin{align} \label{opt: Reg SDP}
															&\min_{X \in \S_n^+} \Tr(\tilde{C}^\top X) + \frac{\rho}{2}  \left| \left|X \right| \right|_2^2 \\ \nonumber
															\text{subject to: } & \Tr(A_i^\top X)=b_i \text{ for } i \in \{1,\dots,m\},
														\end{align}
													\text{} \vspace{-0.4cm}
													
													\noindent where $\tilde{C}=-C \in \S_n$, which is a standard SDP with Tikhonov regularization.

													\ifthenelse{\boolean{longver}}{
															\begin{table*}[t]
																	\centering
																	 \scalebox{0.5}{
																		 		{\begin{tabular}{@{}
																					 				l            
																					 				l            
																					 				r            
																					 				r            
																					 				r            
																					 				r            
																					 				r            
																					 				r            
																					 				r            
																					 				r            
																					 				@{}}
																				 			\toprule
																				 			SDP problem & Metric 
																				 			& Mosek & Alg.~\ref{alg:Dual Gradient Descent} \texttt{scal}=-1 
																				 			& \multicolumn{3}{c}{Alg.~\ref{alg:Dual Gradient Descent} \texttt{scal}=0} 
																				 			& \multicolumn{3}{c}{Alg.~\ref{alg:Dual Gradient Descent} \texttt{scal}=1} \\
																				 			\cmidrule(lr){3-3}\cmidrule(lr){4-4}
																				 			\cmidrule(lr){5-7}\cmidrule(lr){8-10}
																				 			& 
																				 			&       &     
																				 			& $k=[0.2n]$ & $k=[0.1n]$ & $k=[0.05n]$
																				 			& $k=[0.2n]$ & $k=[0.1n]$ & $k=[0.05n]$ \\
																				 			\midrule
																				 			
																				 			\multirow{3}{*}{\shortstack{Degree 4\\{\scriptsize $[n,m]=$[3026, 715]}}}
																				 			& Time (s) 
																				 			&   0.176 &   1.65 
																				 			&    1.8   &    1.5   &   1.46  
																				 			&    2.1   &    1.8   &    1.6 \\
																				 			& $\displaystyle\sqrt{\sum_i(\mathrm{Tr}(A_i^T X)-b_i)^2}$
																				 			& $1.75\times10^{-9}$ & $4.5\times10^{-3}$
																				 			&   2.08  &   5.50  &   2.08  
																				 			&   0.10  &   0.15  &   0.14 \\
																				 			& $|\gamma - \gamma^*|$
																				 			& $6.52\times10^{-10}$ & $1.54\times10^{-5}$
																				 			& $3.1\times10^{-5}$ & $9.56\times10^{-4}$ & $1.0\times10^{-3}$ 
																				 			& $4.5\times10^{-5}$ & $1.76\times10^{-5}$ & $5.1\times10^{-5}$ \\
																				 			\addlinespace
																				 			
																				 			\multirow{3}{*}{\shortstack{Degree 6\\{\scriptsize $[n,m]=$[48401, 5005]}}}
																				 			& Time (s) 
																				 			&   8.82   &  19.26 
																				 			&  16.2    &  11.9    &  10.7   
																				 			&  22.1    &  17.8    &  16.9 \\
																				 			& $\displaystyle\sqrt{\sum_i(\mathrm{Tr}(A_i^T X)-b_i)^2}$
																				 			& $7\times10^{-12}$ & $4.5\times10^{-2}$
																				 			&  2.37   &  2.37    &  2.33   
																				 			&  0.33   &  0.44    &  0.43 \\
																				 			& $|\gamma - \gamma^*|$
																				 			& $1.9\times10^{-12}$ & $2.7\times10^{-4}$
																				 			& $3.53\times10^{-4}$ & $2.10\times10^{-3}$ & $9.9\times10^{-3}$ 
																				 			& $1.25\times10^{-4}$ & $2.3\times10^{-4}$ & $1.43\times10^{-4}$ \\
																				 			\addlinespace
																				 			
																				 			\multirow{3}{*}{\shortstack{Degree 8\\{\scriptsize $[n,m]=$[511226, 24310]}}}
																				 			& Time (s) 
																				 			& 690.6    & 201.0  
																				 			& 160.8    & 134.7    & 114.6   
																				 			& 188.0    & 156.8    & 141.7 \\
																				 			& $\displaystyle\sqrt{\sum_i(\mathrm{Tr}(A_i^T X)-b_i)^2}$
																				 			& $6.83\times10^{-11}$ & 0.17
																				 			&  3.73   &  9.78    &  9.75   
																				 			&  3.37   &  3.10    &  3.30 \\
																				 			& $|\gamma - \gamma^*|$
																				 			& $1.1\times10^{-11}$ & $1.7\times10^{-3}$
																				 			& $1.07\times10^{-4}$ & $6.4\times10^{-3}$ & $1.9\times10^{-2}$ 
																				 			& $1.7\times10^{-3}$ & $1.6\times10^{-3}$ & $1.2\times10^{-3}$ \\
																				 			\addlinespace
																				 			
																				 			\multirow{3}{*}{\shortstack{Degree 10\\{\scriptsize $[n,m]=$	[4008005, 92378]}}}
																				 			& Time (s) 
																				 			& $\infty$ & 2275.7
																				 			& 2179.0   & 1909.5   & 1734.6  
																				 			& 2389.9   & 1992.8   & 1820.6 \\
																				 			& $\displaystyle\sqrt{\sum_i(\mathrm{Tr}(A_i^T X)-b_i)^2}$
																				 			& $\infty$ & 3.5e-2
																				 			&  1.05    &  5.27    &  1.05   
																				 			&  0.86    &  0.99    &  1.20 \\
																				 			& $|\gamma - \gamma^*|$
																				 			& $\infty$ & 4.2e-5
																				 			&  4.8e-8  &  3.0e-4  &  1.3e-3 
																				 			&  1.19e-4 &  7.5e-5  &  1.19e-4 \\
																				 			\bottomrule
																				 		\end{tabular}}
																		}
																	\caption{ \footnotesize Performance of various algorithms for solving the SDP Problem~\eqref{opt: Reg SDP}, {with decision variable size $n \in \N$ and number of constraints $m \in \N$,} associated with Opt.~\eqref{opt: global poly min}. Values of $\infty$ signify the solver failed to provide a solution due to excessive memory requirements.  } 
																	\label{table: bound polynomial} \end{table*}}

													\vspace{-0.2cm}
													\subsection{{SDP} Numerical Examples}
													\ifthenelse{\boolean{longver}}{}{
															\vspace{-0.2cm}}
													The following numerical examples demonstrate the properties of Alg.~\ref{alg:Dual Gradient Descent}. All computation was carried out using an Apple M1 Macbook Pro with 16GB of RAM. 	
													While our numerical experiments include Mosek as a baseline reference, owing to its widespread use as a second-order SDP solver and its integration with SOSTOOLS, we stress that the primary contribution of this work is the development and analysis of RNLA-based approximate PSD projections, rather than proposing an alternative SDP solver. We believe there is room to improve the computation times of the following numerical experiments as we do not take advantage of sparse matrix multiplication, LAPACK implementation or parallel computing. 
													\ifthenelse{\boolean{longver}}{}{
															\vspace{-0.2cm}} \begin{ex} \label{ex: 2d grad descent plot}
															{This example illustrates the impact of using approximate versus exact PSD projections on the convergence behaviour of gradient descent applied to a simple SDLS Problem~\eqref{opt: SDLS} . The problem is designed with only two constraints ($m=2$), allowing the dual objective landscape and gradient descent iterations to be visualized in a two dimensional plot (Figure \ref{fig: 2d grad descent}).}
															
															Consider an SDLS Problem~\eqref{opt: SDLS} with two constraints ($m=2$). We formulate the problem using Matlab to generate $\tilde{A}_i=\texttt{randn(n)}$ and making this normalized and symmetric by setting $A_i=\frac{\tilde{A}_i^\top + \tilde{A_i}}{2 \|\tilde{A}_i\|_2}$ for $i \in \{1,2\}$. To ensure the problem is feasible we randomly generate $X_0 \in \S_n$ 
												   and set $b_i= \Tr(A_i^\top X_0)$. Similarly, we randomly generate $C \in \S_n$. 
												   {Figure~\ref{fig: 2d grad descent} shows the outcome of applying Alg.~\ref{alg:Dual Gradient Descent} using both the full deterministic eigen-decomposition (red stars) and a low-rank randomized approximation (blue diamonds).} Specifically, the blue diamonds show $100$ GD iterations of Alg.~\ref{alg:Dual Gradient Descent}, randomly initiated at $y=[ 1.2753,    0.2418]^\top$ and executed for {$n=200$}, $\rho=1$, $M=100$, $\beta=0.5$, {$k=75$}, $l=10$, {$q=6$}, {scal$=1$}. The red stars show the output of Alg.~\ref{alg:Dual Gradient Descent} when the full deterministic eigen-decomposition is executed during each iteration. We also plot the level sets of the dual objective function $\theta(y):= L(X^*(y),y)$. The algorithm using the full eigen-decomposition converges to the optimal solution, $y^*={[-0.04   -1.19]^\top}$, while the algorithm using low-rank randomized eigenvalue decompositions converges {within a neighbourhood of the true solution}. {Increasing the approximate rank, $k$, tightens the neighbourhood at the cost of increased per-iteration computational expense.} Upon termination, the constraint error, $\sqrt{\sum_{i=1}^2 (\Tr(A_i^\top X) - b_i)^2}$, was { 2.8e-15} for the deterministic GD and {0.0437} for the RNLA GD output. 
															
															\begin{figure}
																	\centering
																	{
																			\includegraphics[width=0.6\linewidth, trim = {1cm 1cm 0.9cm 1cm}, clip]{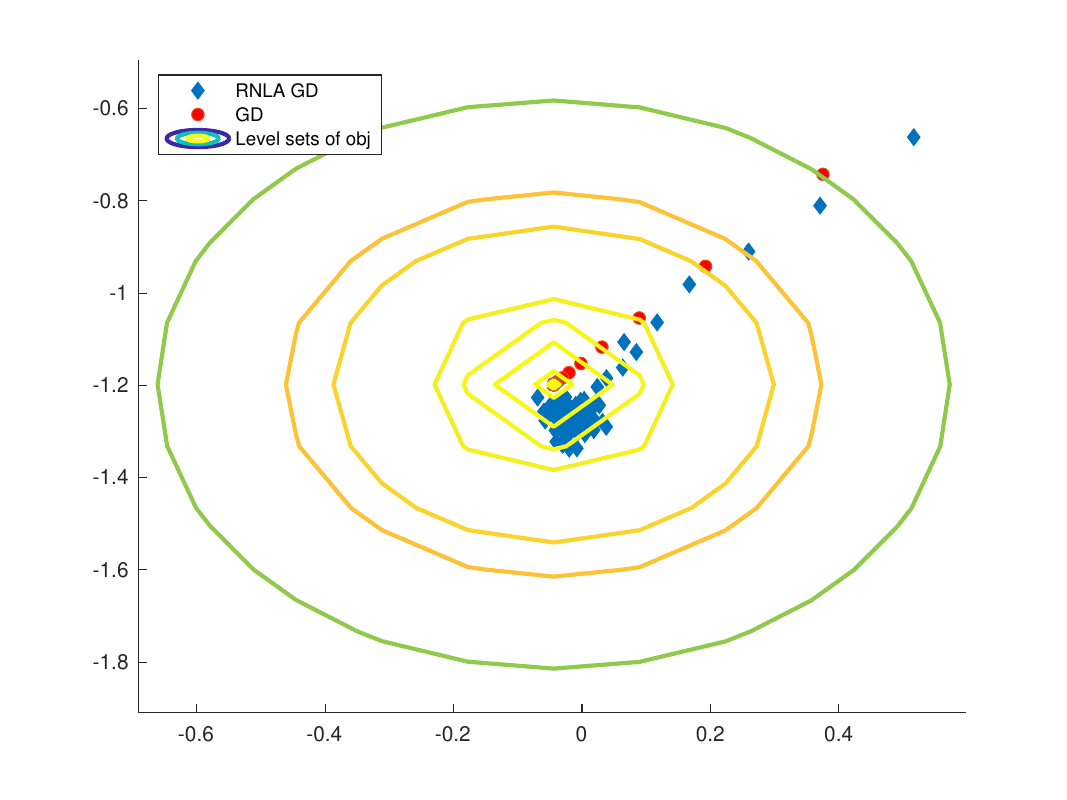}
																			
																		}
																	\vspace{-10pt}
																	\caption{\footnotesize Plot associated with Ex.~\ref{ex: 2d grad descent plot} showing convergence of GD algorithms.} \label{fig: 2d grad descent} 
																\end{figure}
														\end{ex}
													\ifthenelse{\boolean{longver}}{}{
															\vspace{-0.35cm}}
													\begin{ex}
															{We now consider the Sum-of-Squares (SOS) relaxation to an unconstrained polynomial optimization problem~\cite{parrilo03}. The objective is to find a certifiable lower bound on the minimum of a  given polynomial $p:\R^n \to \R$. The SOS relaxation of this is:}
															 \begin{align} \label{opt: global poly min}
																	&	\gamma^*=  \min~ \gamma 
																	\quad \mathrm{subject~ to~ } p(x) + \gamma \in \sum_{SOS,{d}},
																\end{align}
															\ifthenelse{\boolean{longver}}{
																}{
																	\text{ } \\
																	\vspace{-0.8cm}
																}
															{Where $\sum_{SOS,{d}}$ denotes the set of {degree $d \in \N$} SOS polynomials. One of the main reasons  SOS problems are convenient is they can be solved via semidefinite programming~\cite{parrilo03}. Unfortunately the size of the resulting SDP scales unfavorably with the degree of $p$.}
															{We  solved a regularized} version of this SDP, given in Eq.~\eqref{opt: Reg SDP}, using Alg.~\ref{alg:Dual Gradient Descent}. To test performance we randomly generate $p$ such that the solution to Opt.~\eqref{opt: global poly min} is analytically known. We set $p(x)=z_d(x)Xz_d(x) -\pi$, where $z_d$ is the vector of $d$-degree monomial basis functions and $X \in \S_n^+$ is randomly generated. The  solution is then $\gamma^*=\pi$.
															
															The {number of independent variables in } $p$ is fixed to $9$ and the degree is varied. Columns $k=0.2n$, $k=0.1n$, and $k=0.05n$ in Table~\ref{table: bound polynomial} show information about Alg.~\ref{alg:Dual Gradient Descent} when the parameters parameters are fixed to $l=10$, $q=8$, $M=4000$, $N=10$, $\rho=0.1$, $\beta=0.01$, $\eps=1e-10$ and scal={$-1$ or} $1$ or $0$.  {Column Mosek shows the result of applying MOSEK~\cite{aps2019mosek} to the SDP.} Column {scal=$-1$} shows the results of executing Alg.~\ref{alg:Dual Gradient Descent} when the {exact} analytical PSD projection is executed during each iteration of GD{, whereas the remaining six columns show results of executing Alg.~\ref{alg:Dual Gradient Descent} when our inexact approximate PSD projection is used.} 
															
															
															
															{As expected, Mosek, a second-order interior point solver, achieves faster convergence and higher accuracy when the problem size is small enough to fit within available memory. However, as the problem size increases (degree 10 and above), Mosek fails due to its large per-iteration memory requirements, while the first-order gradient-based method remains able to compute solutions.
																	
																When using this first-order method, we observe that replacing the exact PSD projection step with our proposed RNLA-based inexact projection leads to solutions of similar relative quality, both in terms of objective value and constraint satisfaction. This suggests that the RNLA-based projection offers a practical trade-off, maintaining solution accuracy while reducing the per-iteration computational cost when integrated into first order SDP solvers. }

														\end{ex}
														
															In future work, we plan to analyze the end-to-end convergence effects of integrating our approximate PSD projections into various first-order solvers\footnote{ {We have already analysed the end to end convergence of an accelerated deterministic version of Alg.~\ref{alg:Dual Gradient Descent} in~\cite{machado2026lifting}.}}. Additionally, we will investigate adaptive parameter selection strategies for our RNLA algorithms. By dynamically adjusting these parameters during optimization iterations, we aim to simultaneously enhance projection accuracy and overall convergence rates. Finally, we will investigate implementations harnessing sparse matrix multiplications and GPU accelerated architecture.
				
				\fi

%
%
%
\end{document}